\documentclass[3pt]{article}

\usepackage{graphicx}\usepackage{epstopdf}\usepackage{color}
\usepackage{amsmath}\usepackage{amsfonts}\usepackage{amssymb}
 \usepackage{amsthm}
\newtheorem{theorem}{Theorem}[section]
\newtheorem{lemma}{Lemma}[section]

\numberwithin{equation}{section}

\begin{document}

\title{The adaptive Crouzeix-Raviart element method for convection-diffusion eigenvalue problems}

\author{ {Yingyu Du, Qinghua Chen\thanks{Corresponding Author}} \\\\%
{\small School of Mathematical Sciences, }\\{\small
Guizhou Normal University,  Guiyang,  $550001$,  China}\\{\small
1034062167@qq.com, cqh040522@163.com}}
\date{~}
\pagestyle{plain} \textwidth 145mm \textheight 215mm \topmargin 0pt
\maketitle

\indent{\bf\small Abstract~:~} {\small The convection-diffusion eigenvalue problems are hot topics, and computational mathematics community and physics community are concerned about them in recent years. In this paper, we consider the a posteriori error analysis and the adaptive algorithm of the Crouzeix-Raviart nonconforming element method for the convection-diffusion eigenvalue problems. We give the corresponding a posteriori error estimators, and prove their reliability and efficiency. Finally, the numerical results validate the theoretical analysis and show that the algorithm presented in this paper is efficient.}\\
\indent{\bf\small Keywords~:~\scriptsize} {\small
convection-diffusion  eigenvalue problems, the Crouzeix-Raviart element, a posteriori error analysis, adaptive algorithm }
\section{Introduction}
\indent The convection-diffusion eigenvalue problems have a strong background in physics, such as the distribution of contaminated material in nuclear waste pollution. Thus, using finite element methods to solve convection-diffusion eigenvalue problems has attracted much attention of scholars. \cite{J.G4,R.R6,YLHH7} discussed a posteriori error estimates and the adaptive algorithms, \cite{C.Carstensen4} an adaptive homotopy approach, \cite{TLYF9,YHS10} extrapolation methods, \cite{A.N} function value recovery algorithms, \cite{JHYY5} spectral element methods, \cite{YH,pzlan} multilevel correction method, and so on. This paper aims at deriving the a posteriori error estimators and the adaptive algorithm of the Crouzeix-Raviart element(C-R element) methods for the convection-diffusion eigenvalue problems.\\
\indent The adaptive finite element method is a mainstream in scientific computing (see
\cite{ainsworth1,babuska3,verfurth1,shi2}). In past years, the research of the a posteriori error and the adaptive algorithm of
convection-diffusion eigenvalue problems used to adopt the conforming finite element methods(see\cite{YLHH7,JHYY5,babuska3}). \cite{C. Carstensen2} and \cite{Li} discussed a posteriori error estimate of the nonconforming methods for Laplace equation and Laplace eigenvalue problem, respectively. Based on the study of \cite{C. Carstensen2,Li}, this paper first discusses the nonconforming finite element adaptive method for convection-diffusion eigenvalue problems.
  We give the a posteriori error estimators and  prove their reliability and efficiency, and give the adaptive algorithm.
Finally we use some numerical examples to verify our theoretical results.\\
\indent In this paper, $C$ is a positive constant independent of $h$, which may not be the same constant in different places. For simplicity, we use symbol $a\lesssim b$ to replace $a \leq Cb$. The notation $a\approx b$ abbreviates $a\lesssim b\lesssim a$.

\section{Preliminaries}

Consider the following convection-diffusion eigenvalue problem:
\begin{eqnarray}\label{s2.1}
-\Delta u+\mathbf{b}\cdot\nabla u=\lambda u ~~~in~\Omega, ~~~~~ u=0 ~~~ on~\partial \Omega,
\end{eqnarray}
where $\Omega\subset\mathbb{R}^{2}$ is a polygon bounded domain with
boundary $\partial\Omega$.\\
\indent Let
\begin{eqnarray}\label{s2.2}
a(u,v)=\int_\Omega \nabla u\cdot \nabla \overline{v}+\mathbf{b}\cdot \nabla u\overline{v} ~dx,~~~~~b(u,v)=\int_\Omega  u\overline{v} dx.
\end{eqnarray}
The variational problem associated with (\ref{s2.1}) is given by: Find $(\lambda,u)\in\mathbb{C}\times H_{0}^{1}(\Omega)$, $\parallel
u \parallel_{L^{2}(\Omega)}=1,$ such that
\begin{eqnarray}\label{s2.3}
&&a(u,v)=\lambda b(u,v),~~~\forall v\in H_{0}^{1}(\Omega).
\end{eqnarray}
\indent Let $\mathcal{T}_h=\{K\}$ be a regular triangular mesh of
$\Omega$.\\
\indent Let $V_h$ denote the Crouzeix-Raviart nonconforming finite element space over $\mathcal{T}_h$. Then, the C-R element approximation of (\ref{s2.3}) is given as follows: Find $(\lambda_h,u_h)\in \mathbb{C}\times V_h$, $\parallel u_{h} \parallel_{L^{2}(\Omega)}=1$, such that
\begin{eqnarray}\label{s2.4}
a_h(u_h,v)=\lambda_h b(u_h,v),~~~\forall v\in V_h.
\end{eqnarray}
where
\begin{eqnarray}\label{s2.5}
a_h(u_h,v)=\sum\limits_{K}\int_K \nabla_h u_h\nabla \overline{v} +\mathbf{b}\cdot \nabla_h u_h \overline{v} dx.
\end{eqnarray}
Since the discrete space $V_h$ is nonconforming, we regard $\nabla_h$ as the gradient operator which is defined elementwise.\\
\indent The dual problem of (\ref{s2.1}) is as below:
\begin{eqnarray}\label{s2.6}
-\Delta u^*-\nabla\cdot(\overline{\mathbf{b}}u^*)=\lambda^* u^* ~~~in~\Omega, ~~~~~ u^*=0 ~~~ on~\partial \Omega.
\end{eqnarray}
The corresponding variational form of (\ref{s2.6}) is as follows: Find $(\lambda^*,u^*)\in\mathbb{C}\times H_{0}^{1}(\Omega)$, $\parallel u^{*} \parallel_{L^{2}(\Omega)}=1$, such that
\begin{eqnarray}\label{s2.7}
&&a(v,u^*)=\overline{\lambda^*} b(v,u^*),~~~\forall v\in H_{0}^{1}(\Omega),\\\nonumber
\end{eqnarray}
where
\begin{eqnarray}\label{s2.8}
a(v,u^*)=\int_\Omega\nabla v \cdot{\nabla \overline {u^*}}+\nabla v\cdot {\mathbf{b} \overline {u^*}}dx,~~~b(v,u^*)=\int_\Omega  v\overline{u^*} dx.\\\nonumber
\end{eqnarray}
Then the C-R element approximation of (\ref{s2.7}) is as below: Find $(\lambda_h^*,u_h^*)\in \mathbb{C}\times V_h$, $\parallel u_h^* \parallel_{L^{2}(\Omega)}=1$,  such that
\begin{eqnarray}\label{s2.9}
a_h(v,u_h^*,)=\overline{\lambda_h^*} b(v,u_h^*),~~~\forall v\in V_h,
\end{eqnarray}
where
\begin{eqnarray}\label{s2.10}
a_h(v,u_h^*)=\sum\limits_{T}\int_{T}\nabla v \nabla\overline{ u_h^*}+\nabla v\cdot \mathbf{b}\overline{ u_h^*}dx,~~~b(v,u_h^*)=\int_\Omega  v\overline{u_h^*}dx.
\end{eqnarray}

\indent \cite{YHB} discusses the non-conforming finite element
approximation, and proves the error estimates of the discrete
eigenvalues obtained by the Adini element, Morley-Zienkiewicz
element et. al.
 Due to the reference \cite{YHB}, we can deduce the following Lemma.
\begin{lemma}
 For the C-R nonconforming finite element methods of problem (\ref{s2.1}) and (\ref{s2.6}), the a priori error estimates are given:
\begin{eqnarray}\label{s2.11}
&&\parallel\nabla_h(u-u_h)\parallel_{L^2(\Omega)}\lesssim h^{r},\\
\label{s2.12}
&&\parallel u-u_h\parallel_{L^2(\Omega)}\lesssim h^{r}\parallel\nabla_h(u-u_h)\parallel_{L^2(\Omega)},\\
\label{s2.13}
&&\parallel\nabla_h(u^*-u_h^*)\parallel_{L^2(\Omega)}\lesssim h^{\frac{r}{\alpha}},\\
\label{s2.14}
&&\parallel u^{*}-u_h^{*}\parallel_{L^2(\Omega)}\lesssim (h^r\parallel\nabla_h(u^{*}-u_h^{*})\parallel_{L^2(\Omega)})^\frac{1}{\alpha},\\
\label{s2.15}
&&\mid\lambda-\lambda_h\mid\lesssim\parallel\nabla_h(u-u_h)\parallel_{L^2(\Omega)}\cdot\parallel\nabla_h(u^*-u_h^*)\parallel_{L^2(\Omega)}.
\end{eqnarray}
\end{lemma}
Owing to the above conclusions, we can get the following estimate:
there exist some positive constants $0<\beta <1$ and $h_0>0$ (when $h<h_0$) with
\begin{eqnarray}\label{s2.16}\nonumber
&&\mid\lambda-\lambda_h\mid\parallel u\parallel_{L^2(\Omega)}+\mid\lambda_h\mid\parallel u-u_h\parallel_{L^2(\Omega)}+\parallel u-u_h\parallel_{L^2(\Omega)}\\
&&~~~\leq \beta\parallel\nabla_h(u-u_h)\parallel_{L^2(\Omega)}.
\end{eqnarray}

\section{A posteriori error analysis }
\indent Now we introduce some symbols for reading convenience. Suppose $K$ is one given element of $\mathcal{T}_h$, and $h_K$ represents the diameter of $K$. We use $\varepsilon$ to denote the set of all edges in $\mathcal{T}_h$, $\varepsilon(\Omega)$ the set of interior edges and $\varepsilon(K)$ the set of edges of the element $K$, respectively. For any given edge $E\in\mathcal{\varepsilon}(\Omega)$ with length $h_E=|E|$, we assign the fixed unit normal $\nu_E:=(\nu_1,\nu_2)$ and
tangential vector $\tau_E:=(-\nu_2,\nu_1)$. Once $\nu_E$ and $\tau_E$ have been fixed on $E$, in relation to $\nu_E$ one defines the elements $K_{-}\in\mathcal{T}_h$ and $K_{+}\in\mathcal{T}_h$, with $E=K_{+}\bigcap K_{-}$ and $\omega_E=K_{+}\bigcup K_{-}$. Given $E\in\varepsilon(\Omega)$, we denote by $[v]:=(v|_{K_{+}})|_E-(v|_{K_{-}})|_E$ the jump of some $R^{d}$-valued function $v$ defined in $\Omega$  across $E$ with $d=1,2$. And throughout this paper,
$[\cdot]$ denotes the jump of the piecewise smooth function across the
internal edge $E$, and the trace for the boundary edge $E$.\\
 \indent Define the a posteriori error estimators on the element $K$ as below:
 \begin{eqnarray*}
&&\eta_{h,K}:=(h_K^2\parallel\lambda_h u_h+\Delta_h u_h-\mathbf{b}\cdot\nabla_h u_h\parallel_{L^2(K)}^2)^\frac{1}{2},\\
&&\eta_{h,K}^*:=(h_K^2\parallel\lambda_h^* u_h^*+\Delta_h u_h^*+\nabla_h\cdot\mathbf{b} u_h^*\parallel_{L^2(K)}^2)^\frac{1}{2},\\
&&\eta_{h,K,\nu_E}:=(\frac{1}{2}\sum_{E\in\partial K}h_E\parallel[\nabla_h u_h]\cdot\nu_E\parallel_{L^2(E)}^2)^\frac{1}{2},\\
&&\eta_{h,K,\nu_E}^*:=(\frac{1}{2}\sum_{E\in\partial K}h_E\parallel[\nabla_h u_h^*+\mathbf{b} u_h^*]\cdot\nu_E\parallel_{L^2(E)}^2)^\frac{1}{2},\\
&&\eta_{h,K,\tau_E}:=(\frac{1}{2}\sum_{E\in\partial K}h_E\parallel[\nabla_h u_h]\cdot\tau_E\parallel_{L^2(E)}^2)^\frac{1}{2},\\
&&\eta_{h,K,\tau_E}^*:=(\frac{1}{2}\sum_{E\in\partial
K}h_E\parallel[\nabla_h u_h^*+\mathbf{b}
u_h^*]\cdot\tau_E\parallel_{L^2(E)}^2)^\frac{1}{2},
\end{eqnarray*}
and the residual sum on $K$ are given by
\begin{eqnarray}\label{s3.1}
&&\eta_h(K)^{2}:=\eta_{h,K}^2+\sum_{E\in\varepsilon(K),E\not\subset
\partial\Omega}\eta_{h,K,\nu_E}^2+\sum_{E\in\varepsilon(K)}\eta_{h,K,\tau_E}^2,\\\label{s3.2}
&&\eta_h^*(K)^{2}:=(\eta_{h,K}^{*})^2+\sum_{E\in\varepsilon(K),
E\not\subset
\partial\Omega
}(\eta_{h,K,\nu_E}^*)^2
+\sum_{E\in\varepsilon(K)}(\eta_{h,K,\tau_E}^*)^2.
\end{eqnarray}
For any $\mathcal{M}_{h}\subset\mathcal{T}_{h}$, define the estimators over $\mathcal{M}_{h}$ by
\begin{eqnarray}\label{s3.3}
\eta_h(\mathcal{M}_{h})^{2}:= \sum_{K\in
\mathcal{M}_{h}}\eta_h(K)^{2},~~~\eta_h^{*}(\mathcal{M}_{h})^{2}:=
\sum_{K\in \mathcal{M}_{h}}\eta_h^{*}(K)^{2}.
\end{eqnarray}

\indent The left parts of this section aim at proving the reliability and the
efficiency of the estimators $\eta_h(\mathcal{T}_{h})$ and
$\eta_h^{*}(\mathcal{T}_{h})$.
 The reliability of the estimators are based on the following lemma (see\cite{shi2,Li}).

\begin{lemma}
 Under the assumption (\ref{s2.16}) there holds 
\begin{eqnarray}\label{s3.4}\nonumber
&&|a_h(u-u_h,u-u_h )|\lesssim\min_{v\in H_0^1(\Omega)}\parallel\nabla_h(u_h-v)\parallel_{L^2(\Omega)}^2\label{s2.11}\\
&&~~~~~~+\sup_{w\in H_0^1(\Omega)}\frac{|b(\lambda_h
u_h,w)-a_h(u_h,w)|}{\parallel
w\parallel_{L^2(\Omega)}}\parallel\nabla(u-v)\parallel_{L^2(\Omega)},\label{s2.12}
\end{eqnarray}
where $(\lambda,u)\in \mathbb{C}\times H_{0}^{1}(\Omega)$ and $(\lambda_h,u_h)\in \mathbb{C}\times V_h$ are the solutions
 to problems(\ref{s2.3})and(\ref{s2.4}), respectively. For the dual problem, it is similar:
\begin{eqnarray}\label{s3.5}\nonumber
&&|a_h(u^*-u_h^*,u^*-u_h^* )|\lesssim\min_{v\in H_0^1(\Omega)}\parallel\nabla_h(u_h^*-v)\parallel_{L^2(\Omega)}^2\label{s2.11}\\
&&~~~~~~+\sup_{w\in H_0^1(\Omega)}\frac{|b(\lambda_h^*
u_h^*,w)-a_h(u_h^*,w)|}{\parallel
w\parallel_{L^2(\Omega)}}\parallel\nabla(u^*-v)\parallel_{L^2(\Omega)}.
\label{s2.12}
\end{eqnarray}
\end{lemma}
\indent{\bf Proof.}
\indent For any $v\in H_{0}^{1}(\Omega)$,
\begin{eqnarray}\label{s3.6}\nonumber
&&|a_h(u-u_h,u-u_h)|=|a_h(u-u_h,u-v+v-u_h)|\\\nonumber
 &&~~~=|a(u,u-v
)-a_h(u_h,u-v)+a_h(u-u_h,v-u_h)|\\\nonumber
 &&~~~=|b(\lambda u-\lambda_h
u+\lambda_h u-\lambda_h u_h,u-u_h+u_h-v)+b(\lambda_h
u_h,u-v)\\\nonumber &&~~+a_h(u-u_h,v-u_h)-a_h(u_h,u-v)|\\\nonumber
&&~~~\leq|b((\lambda-\lambda_h)u,u-u_h)+b(\lambda_h(u-u_h),u-u_h)|\\\nonumber
&&~~~~~~+|b((\lambda-\lambda_h)u,v-u_h)+b(\lambda_h(u-u_h),v-u_h)|\\\nonumber
&&~~~+|a_h(u-u_h,v-u_h)|\\
&&~~~+|b(\lambda_h u_h,u-v)-a_h(u_h,u-v)|.
\end{eqnarray}
\indent Due to (\ref{s2.16}), we can get
\begin{eqnarray}\label{s3.7}\nonumber
&&|b((\lambda-\lambda_h)u,u-u_h)+b(\lambda_h(u-u_h),u-u_h)|\\\nonumber
&&~~~\lesssim(\mid\lambda-\lambda_h\mid\parallel u\parallel_{L^2(\Omega)}+\mid\lambda_h\mid\parallel u-u_h\parallel_{L^2(\Omega)})\parallel u-u_h\parallel_{L^2(\Omega)}\\
&&~~~\leq
\beta^{2}\parallel\nabla_h(u-u_h)\parallel_{L^2(\Omega)}^{2}.
\end{eqnarray}
\begin{eqnarray}\label{s3.8}\nonumber
&&|b((\lambda-\lambda_h)u,v-u_h)+b(\lambda_h(u-u_h),v-u_h)|\\\nonumber
&&~~~\lesssim(\mid\lambda-\lambda_h\mid\parallel u\parallel_{L^2(\Omega)}+\mid\lambda_h\mid\parallel u-u_h\parallel_{L^2(\Omega)})\parallel v-u_h\parallel_{L^2(\Omega)}\\
&&~~~\leq
\beta\parallel\nabla_h(u-u_h)\parallel_{L^2(\Omega)}\parallel
v-u_h\parallel_{L^2(\Omega)}.
\end{eqnarray}
\indent Using the Young and Poincar$\acute{e}$ inequalities, we obtain    
\begin{eqnarray}\label{s3.9}\nonumber
&&~\beta\parallel\nabla_h(u-u_h)\parallel_{L^2(\Omega)}\parallel v-u_h\parallel_{L^2(\Omega)}\\\nonumber
&&~\leq\frac{1}{2}(\varepsilon^2\beta^2\parallel\nabla_h(u-u_h)\parallel_{L^2(\Omega)}^2 +\frac{1}{\varepsilon^{2}}\parallel v-u_{h}\parallel^{2}_{L^{2}(\Omega)})\\
&&~\leq\frac{1}{2}(\varepsilon^2\beta^2\parallel\nabla_h(u-u_h)\parallel_{L^2(\Omega)}^2 +\frac{C1}{\varepsilon^{2}}\parallel \nabla_{h}(v-u_{h})\parallel^{2}_{L^{2}(\Omega)}).
\end{eqnarray}
\indent The inequality  (\ref{s3.6}) gives\\
\begin{eqnarray}\label{s3.10}\nonumber
&&~~~|a_h(u-u_h,v-u_h)|\\\nonumber
 &&=|\sum\limits_{K}\int_K \nabla_h
(u-u_h)\cdot \nabla_h \overline{(v-u_h)}+\mathbf{b}\cdot \nabla_h
(u-u_h)\overline{(v-u_h)}dx|\\\nonumber
&&\lesssim\sum\limits_{K}\{\parallel\nabla_h(u-u_h)\parallel_{L^2(K)}
\parallel\nabla_h(v-u_h)\parallel_{L^2(K)}\\\nonumber
&&~~~+\mid\mathbf{b}\mid\parallel\nabla_h(u-u_h)\parallel_{L^2(K)}
\parallel(v-u_h)\parallel_{L^2(K)}\}\\\nonumber
&&\lesssim\parallel\nabla_h(u-u_h)\parallel_{L^2(\Omega)}\parallel\nabla_h(v-u_h)\parallel_{L^2(\Omega)}\\\nonumber
&&\leq(\frac{1}{2}(\varepsilon^2\parallel\nabla_h(u-u_h)\parallel_{L^2(\Omega)}^2 +\frac{1}{\varepsilon^{2}}\parallel v
-u_{h}\parallel^{2}_{L^{2}(\Omega)})\\
&&\leq(\frac{1}{2}(\varepsilon^2\parallel\nabla_h(u-u_h)\parallel_{L^2(\Omega)}^2
+\frac{C2}{\varepsilon^{2}}\parallel
\nabla_{h}(v-u_{h})\parallel^{2}_{L^{2}(\Omega)}).
\end{eqnarray}
Combining (\ref{s3.7}), (\ref{s3.8}), (\ref{s3.9}) with (\ref{s3.10}), we obtain from (\ref{s3.6})
\begin{eqnarray}\label{s3.11}\nonumber
&&|a_h(u-u_h,u-u_h)|
\leq(\frac{1}{2}\varepsilon^2\beta^2+\beta^2+\frac{C1}{2}\varepsilon^2)\parallel\nabla_h(u-u_h)\parallel_{L^2(\Omega)}^2
\\\nonumber
&&~~~~~~+(\frac{1}{2\varepsilon^2}+\frac{C2}{2\varepsilon^2})
\parallel\nabla_h(u_h-v)\parallel_{L^2(\Omega)}^2+|b(\lambda_h
u_h,u-v)-a_h(u_h,u-v)|,
\end{eqnarray}
then, we have
\begin{eqnarray}\label{s3.12}\nonumber
&&~~~|a_h(u-u_h,u-u_h)|\\\nonumber
 &&\lesssim
\parallel\nabla_h(u_h-v)\parallel_{L^2(\Omega)}^2+|b(\lambda_h
u_h,u-v)-a_h(u_h,u-v)|\\\nonumber
 &&\lesssim
\parallel\nabla_h(u_h-v)\parallel_{L^2(\Omega)}^2\\\nonumber
&&~~~+\frac{|b(\lambda_h
u_h,u-v)-a_h(u_h,u-v)|}{\parallel\nabla(u-v)\parallel_{L^2(\Omega)}}\parallel\nabla(u-v)\parallel_{L^2(\Omega)}\\\nonumber
&&\lesssim\min_{v\in H_0^1(\Omega)}\parallel\nabla_h(u_h-v)\parallel_{L^2(\Omega)}^2\label{s2.11}\\
&&~~~+\sup_{w\in H_0^1(\Omega)}\frac{|b(\lambda_h
u_h,w)-a_h(u_h,w)|}{\parallel
w\parallel_{L^2(\Omega)}}\parallel\nabla(u-v)\parallel_{L^2(\Omega)}.\label{s2.12}
\end{eqnarray}
Then the proof of (\ref{s3.4}) is finished, and the proof of (\ref{s3.5}) is similar.~~$\Box$

\indent Based on the work of \cite{Li,C. Carstensen8}, we have the following Lemma:
\begin{lemma}
The following estimate is valid:
 \begin{eqnarray}\label{s3.13}
\min_{v\in H_0^1(\Omega)}{\parallel\nabla_h (u_h-v)\parallel}^{2}_{L^{2}(\Omega)} \lesssim\sum_{E\in\varepsilon}h_E\parallel[\nabla_h u_h]\cdot\tau_E\parallel_{L^2(E)}^2.
\end{eqnarray}
\end{lemma}
Let $S_0^1(\mathcal{T}_h)$ denote the elementwise linear conforming finite element space over $\mathcal{T}_h$. For the analysis in the rear, we need the $\mathrm{Cl\acute{e}ment-type}$ interpolation operator $\mathcal{L}:H_0^1(\Omega)\mapsto S_0^1(\mathcal{T}_h)$ with the properties(see\cite{P.C,C.Carstensen3,C.B})
\begin{eqnarray}\label{s3.14}
&&\parallel\nabla\mathcal{L}\varphi\parallel_{L^2(K)}+\parallel h_K^{-1}(\varphi-\mathcal{L}\varphi)\parallel_{L^2(K)}\lesssim\parallel\nabla\varphi\parallel_{L^2(\omega_K)},
\end{eqnarray}
and
\begin{eqnarray}\label{s3.15}
&&\parallel h_E^\frac{-1}{2}(\varphi-\mathcal{L}\varphi)\parallel_{L^2(E)}\lesssim\parallel\nabla\varphi\parallel_{L^2(\omega_K)},
\end{eqnarray}
\indent where $E\in\mathcal{\varepsilon}(K)$ and $\varphi\in H_0^1(\Omega)$. In this paper, $\omega_K$ denotes the element patch defined as
\begin{eqnarray}\label{s3.16}
&&\omega_K:=\{{T\in\mathcal{T}_h:\overline{T}\bigcap\overline{K}\neq\emptyset}\}.\\\nonumber
\end{eqnarray}
Refering to \cite{Li}, we can prove the following Lemma.
\begin{lemma}
The following estimates are valid:
\begin{eqnarray}
\label{s3.17}
\sup_{w\in H_0^1(\Omega)}\frac{|b(\lambda_h
u_h,w)-a_h(u_h,w)|}{\parallel\nabla
w\parallel_{L^2(\Omega)}}\lesssim(\sum\limits_{K\in
\mathcal{T}_h}\eta_{h,K}^2+\sum_{E\in\varepsilon(\Omega)}\eta_{h,K,\nu_E}^2)^{\frac{1}{2}},
\end{eqnarray}
\begin{eqnarray}
\label{s3.18}
\sup_{w\in H_0^1(\Omega)}\frac{|b(w,\lambda_h^*
u_h^*)-a_h(w,u_h^*)|}{\parallel\nabla w\parallel_{L^2(\Omega)}}
\lesssim(\sum\limits_{K\in
\mathcal{T}_h}(\eta_{h,K}^{*})^2+\sum_{E\in\varepsilon(\Omega)}(\eta_{h,K,\nu_E}^*)^2)^{\frac{1}{2}}.
\end{eqnarray}
\end{lemma}
\indent{\bf Proof.}
\indent
Using the estimates (\ref{s3.14}) and (\ref{s3.15}) and integrating by parts, we can deduce that
\indent   
\begin{eqnarray}\label{s3.19}
&&\mid b(\lambda_h u_h,w)-a_h(u_h,w)\mid=|(\lambda_h
u_h,w-\mathcal{L}w)_{L^2(\Omega)}-a_h(u_h,w-\mathcal{L}w)_{L^2(\Omega)}|\nonumber\\
 &&~~~=|\sum_K\int_K\lambda_h
u_h\overline{(w-\mathcal{L}w)}dx-\sum_K\int_K-\Delta_h
u_h\overline{(w-\mathcal{L}w)}\nonumber\\
 &&~~~~~~+\mathbf{b}\cdot\nabla_h
u_h\overline{(w-\mathcal{L}w)}dx-\int_{\partial k }\frac{\partial u_h}{\partial
\nu}\overline{(w-\mathcal{L}w)}ds|\nonumber\\
 &&~~~=|\sum_K\int_K(\lambda_h
u_h+\Delta_h u_h-\mathbf{b}\cdot\nabla_h
u_h)\overline{(w-\mathcal{L}w)}dx\nonumber-\sum_K\int_{\partial k
}\frac{\partial u_h}{\partial \nu}\overline{(w-\mathcal{L}w)}ds|\nonumber\\
&&\lesssim\sum_K\parallel\lambda_h u_h+\Delta_h
u_h-\mathbf{b}\cdot\nabla_h u_h\parallel_{L^2(K)}\cdot
h_K\parallel\nabla w\parallel_{L^2(\omega_K)}\nonumber\\
&&~~~~~~~~~+\sum_{E\in\varepsilon(\Omega)}\parallel\nabla_h
u_h\cdot\nu_E\parallel_{L^2(E)}\cdot h_E^{\frac{1}{2}}\parallel\nabla
w\parallel_{L^2(\omega_K)}\nonumber\\
 &&\lesssim(\sum_K
h_K^2\parallel\lambda_h u_h+\Delta_h u_h-\mathbf{b}\cdot\nabla_h
u_h\parallel_{L^2(K)}^2)^{\frac{1}{2}}(\sum_K\parallel\nabla
w\parallel_{L^2(\omega_K)}^2)^{\frac{1}{2}}\nonumber\\
&&~~~~~~~~~+(\sum_{E\in\varepsilon(\Omega)}h_E\parallel\nabla_h
u_h\cdot\nu_E\parallel_{L^2(E)}^2)^{\frac{1}{2}}(\sum_{E\in\varepsilon(\Omega)}\parallel\nabla
w\parallel_{L^2(\omega_K)}^2)^{\frac{1}{2}}\nonumber\\
&&\lesssim(\sum_K h_K^2\parallel\lambda_h u_h+\Delta_h
u_h-\mathbf{b}\cdot\nabla_h u_h\parallel_{L^2(K)}^2\nonumber\\
&&~~~~~~~~~+\sum_{E\in\varepsilon(\Omega)}h_E\parallel[\nabla_h
u_h]\cdot\nu_E\parallel_{L^2(E)}^2)^{\frac{1}{2}}
\parallel\nabla w\parallel_{L^{2}(\Omega)}\nonumber\\
&&\lesssim(\sum\limits_{K\in
\mathcal{T}_h}\eta_{h,K}^2+\sum_{E\in\varepsilon(\Omega)}\eta_{h,K,\nu_E}^2)^{\frac{1}{2}}\parallel\nabla
w\parallel_{L^{2}(\Omega)}.
\end{eqnarray}
This ends the proof. The proof of (\ref{s3.18}) is similar.~~$\Box$

\indent Combining Lemma 3.2 with Lemma 3.3, we can get the reliability of the a posteriori error estimators.
\begin{theorem}
 Let $(\lambda,u)\in \mathbb{C}\times H_{0}^{1}(\Omega)$ and $(\lambda_h,u_h)\in \mathbb{C}\times V_h$ be the solutions to problems
 (\ref{s2.3}) and (\ref{s2.4}), and let
 $(\lambda^*,u^*)\in \mathbb{C}\times H_{0}^{1}(\Omega)$ and $(\lambda_h^*,u_h^*)\in \mathbb{C}\times V_h$ be the solutions to problems(\ref{s2.7})
 and(\ref{s2.9}), respectively. Under the assumption (\ref{s2.16}) there holds
 \begin{eqnarray}\label{s3.20}
&&\|u-u_h\|_{h}^{2}\lesssim
\eta_h(\mathcal{T}_{h})^2,\\\label{s3.21}
 &&
\|u^*-u_h^*\|_{h}^{2}\lesssim\eta_h^*(\mathcal{T}_{h})^2,\\\label{s3.22}
 &&|\lambda-\lambda_{h}|\lesssim \eta_h(\mathcal{T}_{h})^2+\eta_h^{*}(\mathcal{T}_{h})^2.
\end{eqnarray}
\end{theorem}
\indent{\bf Proof.}
Combining  Lemmas 3.1-3.3 we get (\ref{s3.20}) and
(\ref{s3.21}). Substituting (\ref{s3.20}) and (\ref{s3.21}) into
(\ref{s2.15}) yields (\ref{s3.22}).~~$\Box$

\noindent Next, we shall prove the efficiency of the a posteriori error estimators.
\begin{theorem}
Assume the conditions of Theorem 3.1 hold, then
 \begin{eqnarray}\label{s3.23}
&&\eta_h(\mathcal{T}_{h})^2\lesssim
\|u-u_h\|_{h}^{2},\\\label{s3.24}
&&\eta_h^*(\mathcal{T}_{h})^2\lesssim \|u^*-u_h^*\|_{h}^{2}.
\end{eqnarray}
\end{theorem}
\indent{\bf Proof.}
\indent{\bf 1. Proof of $\sum\limits_{K\in \mathcal{T}_h}\eta_{h,K}^2\lesssim\parallel\nabla(u-u_h)\parallel_{L^2(\Omega)}^2$\\}
\indent~~Given $K\in\mathcal{T}_h$, let $b_K=27\lambda_1\lambda_2\lambda_3$ with $\lambda_i$, $i=1,2,3$. Define
\begin{eqnarray}\label{s3.25}
v_K=b_K(\lambda_h u_h+\Delta_h u_h-\mathbf{b}\cdot\nabla_h u_h)
\end{eqnarray}
Then, we have
\begin{eqnarray}\label{s3.26}\nonumber
&&\parallel\lambda_h u_h+\Delta_h u_h-\mathbf{b}\cdot\nabla_h
u_h\parallel_{L^2(K)}^2\approx(\lambda_h u_h+\Delta_h
u_h-\mathbf{b}\cdot\nabla_h u_h,v_K)_{L^2(K)}\\\nonumber
&&~~~=(\lambda_h u_h-\lambda u+\lambda u+\Delta_h
u_h-\mathbf{b}\cdot\nabla_h u_h,v_K)_{L^2(K)}\\\nonumber
&&~~~=(\lambda_h u_h-\lambda u,v_K)_{L^2(K)}+(-\Delta
u+\mathbf{b}\cdot\nabla u+\Delta_h u_h-\mathbf{b}\cdot\nabla_h
u_h,v_K)_{L^2(K)}\\\nonumber
 &&~~~=(\lambda_h u_h-\lambda
u,v_K)_{L^2(K)}+(-\Delta u+\Delta_h u_h,v_K)_{L^2(K)}\\\nonumber
&&~~~~~~+(\mathbf{b}\cdot\nabla u-\mathbf{b}\cdot\nabla_h
u_h,v_K)_{L^2(K)}\\\nonumber
&&~~~=(\lambda_h u_h-\lambda u,v_K)_{L^2(K)}+(\nabla_h (u-u_h),\nabla v_K)_{L^2(K)}\\
&&~~~~~~+(\mathbf{b}\cdot\nabla u-\mathbf{b}\cdot\nabla_h
u_h,v_K)_{L^2(K)}.
\end{eqnarray}
Using the Young inequalities in (\ref{s3.26}) to obtain
\begin{eqnarray}\label{s3.27}\nonumber
&&|(\lambda_h u_h-\lambda u,v_K)_{L^2(K)}|\leq\parallel\lambda_h
u_h-\lambda u\parallel_{L^2(K)}\parallel
v_K\parallel_{L^2(K)}\\\nonumber
 &&~~~\lesssim\parallel\lambda_h
u_h-\lambda u\parallel_{L^2(K)}\parallel\lambda_h u_h+\Delta
u_h-\mathbf{b}\cdot\nabla_h u_u\parallel_{L^2(K)}\\\nonumber
&&\leq\frac{1}{2}(\frac{1}{\varepsilon^2}\parallel\lambda_h u_h-\lambda u\parallel_{L^2(K)}^2+\varepsilon^2\parallel\lambda_h u_h+\Delta_h u_h-\mathbf{b}\cdot\nabla_h u_u\parallel_{L^2(K)}^2).\\
\end{eqnarray}
Thanks to the assumption (\ref{s3.14}) and using the Young inequalities we can have
\begin{eqnarray}\label{s3.28}\nonumber
&&|(\nabla_h (u-u_h),\nabla v_K)_{L^2(K)}|\leq\parallel\nabla_h (u-u_h)\parallel_{L^2(K)}\parallel\nabla
v_K\parallel_{L^2(K)}\\\nonumber
 &&~~~\lesssim h_K^{-1}\parallel\nabla_h (u-u_h)\parallel_{L^2(K)}\parallel v_K\parallel_{L^2(K)}\\\nonumber
&&~~~\leq\frac{1}{2}(\frac{1}{\varepsilon^2}\cdot h_K^{-2}\parallel\nabla_h (u-u_h)\parallel_{L^2(K)}^2+\varepsilon^2\parallel\lambda_h u_h+\Delta_h u_h-\mathbf{b}\cdot\nabla_h u_u\parallel_{L^2(K)}^2).\\
\end{eqnarray}
and
\begin{eqnarray}\label{s3.29}\nonumber
&&|(\mathbf{b}\cdot\nabla u-\mathbf{b}\cdot\nabla_h
u_h,v_K)_{L^2(K)}|\leq\parallel \mathbf{b}\cdot\nabla
u-\mathbf{b}\cdot\nabla_h u_h\parallel_{L^2(K)}\parallel
v_K\parallel_{L^2(K)}\\\nonumber
&&~~~\leq\frac{1}{2}(\frac{1}{\varepsilon^2}\parallel \mathbf{b}\cdot\nabla u-\mathbf{b}\cdot\nabla_h u_h\parallel_{L^2(K)}^2+\varepsilon^2\parallel\lambda_h u_h+\Delta_h u_h-\mathbf{b}\cdot\nabla_h u_u\parallel_{L^2(K)}^2).\\
\end{eqnarray}
then combining (\ref{s3.27})-(\ref{s3.29}) can yield:
\begin{eqnarray}\label{s3.30}\nonumber
&&\eta_{h,K}^2=h_K^2\parallel\lambda_h u_h+\Delta_h
u_h-\mathbf{b}\cdot\nabla_h u_h\parallel_{L^2(K)}^2\\\nonumber
&&~~~\lesssim\parallel\nabla_h (u-u_h)\parallel_{L^2(K)}^2+h_K^2\parallel\lambda_h
u_h-\lambda u\parallel_{L^2(K)}^2\\\nonumber
&&~~~~~~+h_K^2\parallel \mathbf{b}\cdot\nabla u-\mathbf{b}\cdot\nabla_h u_h\parallel_{L^2(K)}^2. \\
\end{eqnarray}
Then,we have
\begin{eqnarray}\label{s3.31}\nonumber
&&\sum_{K\in\mathcal{T}_h} \eta_{h,K}^2\lesssim\parallel\nabla_h (u-u_h)\parallel_{L^2(\Omega)}^2+\sum_{K\in\mathcal{T}_h}h_K^2\parallel\lambda_h u_h-\lambda u\parallel_{L^2(K)}^2\\\nonumber
&&~~~~~~+\sum_{K\in\mathcal{T}_h}h_K^2\parallel \mathbf{b}\cdot\nabla u-\mathbf{b}\cdot\nabla_h u_h\parallel_{L^2(K)}^2\\
&&~~~\lesssim\parallel\nabla_h(u-u_h)\parallel_{L^2(\Omega)}^2.
\end{eqnarray}
\indent  {\bf~2. Proof of $ \sum\limits_{E\in\varepsilon(\Omega)} \eta_{h,K,\nu_E}^2\lesssim\parallel\nabla_h (u-u_h)\parallel_{L^2(\Omega)}^2$\\}
\indent Given any edge $E\in\varepsilon(\Omega)$, let $b_E\in H_0^1(\omega_E)$ denote the piecewise polynomial function vanishing at the midside point of $E$ \cite{R.V}. Define
\begin{eqnarray}\label{s3.32}
v_E=b_E[\nabla_h u_h]\cdot\nu_E.
\end{eqnarray}
Then we have
\begin{eqnarray}\label{s3.33}\nonumber
&&\parallel[\nabla_h u_h]\cdot \nu_E \parallel_{L^2(E)}^2\approx([\nabla_h u_h]\cdot \nu_E,v_E)_{L^2(E)}\\
&&~~~=\int_{\omega_E}\Delta_h u_h\cdot \overline{v_E}dx+\int_{\omega_E}\nabla_h
u_h\cdot\nabla \overline{v_E}dx.
\end{eqnarray}
Due to
\begin{eqnarray}\label{s3.34}\nonumber
&&\int_{\omega_E}\lambda u \overline{v}dx =\int_{\omega_E}\nabla u\cdot\nabla \overline{v}dx+\int_{\omega_E}\mathbf{b}\cdot\nabla u \overline{v}dx.\\\nonumber
\end{eqnarray}
and (\ref{s3.14}), (\ref{s3.33}) can be estimated as
\begin{eqnarray}\label{s3.35}\nonumber
&&\int_{\omega_E}\nabla_h(u_h-u)\cdot\nabla \overline{v_E}dx-\int_{\omega_E}\mathbf{b}\cdot\nabla u \overline{v_E}dx+\int_{\omega_E}(\lambda u +\Delta_h u_h)\overline{v_E}dx\\\nonumber
&&=\int_{\omega_E}\nabla_h(u_h-u)\cdot\nabla \overline{v_E}dx-\int_{\omega_E}\mathbf{b}\cdot\nabla_h (u-u_h) \overline{v_E}dx\\\nonumber
&&~~~+\int_{\omega_E}(-\mathbf{b}\cdot\nabla_h u_h+\Delta_h u_h+\lambda_h u_h)\overline{v_E}dx
+\int_{\omega_E}(\lambda u -\lambda_h u_h)\overline{v_E}dx\\\nonumber
&&\lesssim\parallel\nabla_h(u_h-u)\parallel_{L^2(\omega_E)}\parallel\nabla v_E\parallel_{L^2(\omega_E)}+\parallel\nabla_h(u_h-u)\parallel_{L^2(\omega_E)}\parallel v_E\parallel_{L^2(\omega_E)}\\\nonumber
&&~~~+\parallel-\mathbf{b}\cdot\nabla_h u_h+\Delta_h u_h+\lambda_h u_h\parallel_{L^2(\omega_E)}\parallel v_E \parallel_{L^2(\omega_E)}\\\nonumber
&&\lesssim h_E^{-1}\parallel\nabla_h(u_h-u)\parallel_{L^2(\omega_E)}\parallel v_E\parallel_{L^2(\omega_E)}\\
&&\lesssim h_E^{-1}\parallel\nabla_h(u_h-u)\parallel_{L^2(\omega_E)}^2.
\end{eqnarray}
Then, we obtain
\begin{eqnarray}\label{s3.36}\nonumber
&&\sum\limits_{E\in\varepsilon(\Omega)} \eta_{h,K,\nu_E}^2=\sum_{E\in\varepsilon(\Omega)}h_E\parallel[\nabla_h u_h]\cdot \nu_E \parallel_{L^2(E)}^2 \\
&&~~~~~~~~~~~~~~~~~~\lesssim \parallel\nabla_h(u_h-u)\parallel_{L^2(\Omega)}^2.
\end{eqnarray}
\indent {\bf~3. Proof of $\sum\limits_{E\in\varepsilon}\eta_{h,K,\tau_E}^2 \lesssim\parallel\nabla_h (u-u_h)\parallel_{L^2(\Omega)}^2$\\}
\indent With the edge bubble function $b_E$ as in (\ref{s3.32}), we define
\begin{eqnarray}\label{s3.37}
v_E=b_E[\nabla_h u_h]\cdot\tau_E.
\end{eqnarray}
Then, we have
\begin{eqnarray}\label{s3.38}\nonumber
&&\parallel[\nabla_h u_h]\cdot \tau_E \parallel_{L^2(E)}^2 \approx([\nabla_h u_h]\cdot \tau_E,v_E)_{L^2(E)}\\
&&~~~~~~~~~~~~~~~~~~~~~~~~~=\int_E[\nabla_h u_h]\cdot\tau_E\cdot \overline{v_E}ds.
\end{eqnarray}
Noting that $\nu_E=(n_x,n_y)$ and $\tau_E=(-n_y,n_x)$, (\ref{s3.38}) can be estimated as
\begin{eqnarray}\label{s3.39}\nonumber
&&\int_E[-(u_h)_x n_y+(u_h)_y n_x]\cdot \overline{v_E}ds\\\nonumber
&&=\int_{\omega_E}-(u_h)_{xy} \overline{v_E}-(u_h)_x (\overline{v_E})_y+(u_h)_{yx}\overline{v_E}+(u_h)_y (\overline{v_E})_xdx\\\nonumber
&&=\int_{\omega_E}\nabla_h(u_h)\cdot curl \overline{v_E}dx  \\
&&=\int_{\omega_E}\nabla_h(u_h-u)\cdot curl \overline{v_E}dx.
\end{eqnarray}
where  $curl~\overline{v_E}=(-(\overline{v_E})_y,(\overline{v_E})_x)$ and $\int_{\omega_E}\nabla u\cdot curl~\overline{v_E}dx=0$.\\
An application of the inverse estimate leads to
\begin{eqnarray}\label{s3.40}\nonumber
&&\sum\limits_{E\in\varepsilon}\eta_{h,K,\tau_E}^2=\sum_{E\in\varepsilon}h_E\parallel[\nabla_h u_h]\cdot \tau_E \parallel_{L^2(E)}^2 \\
&&~~~~~~~~~~~~~~~\lesssim \parallel\nabla_h(u_h-u)\parallel_{L^2(\Omega)}^2.
\end{eqnarray}
Thanks to the following conclusion
\begin{eqnarray}\label{s3.41}
\parallel\nabla_h(u_h-u)\parallel_{L^2(\Omega)}^2\lesssim a_h(u^*-u_h^*,u^*-u_h^*),
\end{eqnarray}
combining (\ref{s3.31}), (\ref{s3.36}) with (\ref{s3.40}), we obtain (\ref{s3.23}). The proof of (\ref{s3.24}) is similar.~~$\Box$

Combining Lemmas $3.1,~3.2,~3.3$ and Theorem $3.2$, we derive the following theorem: 
\begin{theorem}
Let $(\lambda,u)\in\mathbb{C}\times H_0^1(\Omega)$ and $(\lambda_h,u_h)\in\mathbb{C}\times V_h$ be the solution to problems (\ref{s2.3}) and (\ref{s2.4}), respectively. Then
\begin{eqnarray}\label{s3.42}
a_h(u-u_h,u-u_h )\approx \eta_h^2.
\end{eqnarray}
Let $(\lambda^*,u^*)$ and $(\lambda_h^*,u_h^*)$ be the eigenpairs of the adjoint problems (\ref{s2.7}) and (\ref{s2.9}), respectively. Then
\begin{eqnarray}\label{s3.43}
a_h(u^*-u_h^*,u^*-u_h^*)\approx (\eta_h^*)^2.
\end{eqnarray}
\end{theorem}

\section{The adaptive
algorithm and numerical results }

\indent Using the a posteriori error estimates and consulting the
existing standard algorithm (see, e.g., \cite{J.G4,R.R6,YLHH7}),
we obtain the following adaptive algorithm of the C-R element for the convection-diffusion eigenvalue problem (\ref{s2.1}):\\
\indent{\bf Algorithm 1.}\\
\indent Choose parameter $0<\theta<1$.\\
\indent{\bf Step 1.}~Pick any initial mesh $\mathcal{T}_{h_{0}}$ with mesh size $h_{0}$.\\
\indent{\bf Step 2.}~Solve (\ref{s2.4}) and (\ref{s2.9}) on $\mathcal{T}_{h_{0}}$ for discrete solution $(\lambda_{h_{0}}, u_{h_{0}}, u_{h_{0}}^{*})$.\\
\indent{\bf Step 3.}~Let $l=0$.\\
\indent{\bf Step 4.}~Compute the local indicators ${\eta}_{h_{l}}(K)^{2}+{\eta}_{h_{l}}^{*}(K)^{2}$.\\
\indent{\bf Step 5.}~Construct
$\widehat{\mathcal{T}}_{h_{l}}\subset\mathcal{T}_{h_{l}}$ by {\bf
Marking Strategy E}
and parameter $\theta$.\\
\indent{\bf Step 6.}~Refine $\mathcal{T}_{h_{l}}$ to get a new mesh
$\mathcal{T}_{h_{l+1}}$
by Procedure ${\bf Refine}$.\\
\indent{\bf Step 7.}~Solve (\ref{s2.4}) and (\ref{s2.9}) on $\mathcal{T}_{h_{l+1}}$
 for discrete solution $(\lambda_{h_{l+1}}, u_{h_{l+1}}, u_{h_{l+1}}^{*})$.\\
\indent{\bf Step 8.}~Let $l=l+1$ and go to Step 4.\\

\indent {\bf Marking Strategy E}\\
\indent Given parameter $0<\theta<1$:\\
\indent{\bf Step 1.}~~Construct a minimal subset
$\widehat{\mathcal{T}}_{h_{l}}$ of $\mathcal{T}_{h_{l}}$ by
selecting some elements in $\mathcal{T}_{h_{l}}$ such that
\begin{eqnarray*}
\sum\limits_{K\in
\widehat{\mathcal{T}}_{h_{l}}}({\eta}_{h_{l}}(K)^{2}+{\eta}_{h_{l}}^{*}(K)^{2})
\geq
\theta({\eta}_{h_{l}}(\mathcal{T}_{h_{l}})^{2}+{\eta}_{h_{l}}^{*}(\mathcal{T}_{h_{l}})^{2}).
\end{eqnarray*}
\indent{\bf Step 2.}~~Mark all the elements in
$\widehat{\mathcal{T}}_{h_{l}}$.\\

\indent Next, we will present some numerical experiments by using the triangular C-R element . We use MATLAB 2012 together with the package of IFEM \cite{L.chen} to solve the (\ref{s2.4}) and (\ref{s2.9}) as below. For simplicity of the presentation, we use the following notations:\\
\indent $\lambda_{k,h}$: the $k$-th finite element eigenvalue.\\
\indent $\lambda_k$: the $k$-th exact eigenvalue.\\
\indent $\Phi(\lambda_{k,h}):$  the a posteriori error indicator for $\lambda_{k,h}$.\\
\indent $N_{k,l}(b_1):$  number of degrees of freedom for $\lambda_{k,h}$  after the $i$-th iteration when $\mathbf{\mathbf{b}}=(b_1,0)^T$.\\

\indent {\bf Example 1.}~ Let $\Omega=(0,1)^2$ and
$\mathbf{b}=(b_1,b_2)^{T}$. Consider the convection-diffusion
eigenvalue problem (\ref{s2.1}) whose eigenvalues are
\begin{eqnarray*}
\frac{b_1^2+b_2^2}{4}+\pi^2(j^2+i^2),
\end{eqnarray*}
where $j, i\in N_+$. 
 We know that $\lambda_1=\frac{b_1^2+b_2^2}{4}+2\pi^2$, $\lambda_2=\lambda_3=\frac{b_1^2+b_2^2}{4}+5\pi^2$. We restrict our attention to the case of $\mathbf{b}=(1,0)^T, \mathbf{b}=(3,0)^T$, and $\mathbf{b}=(10,0)^T$. Some adaptive refined meshes are shown in Figures 1 and 2 and the numerical results are shown in table 1. From the results we can see that   the a posteriori error indicators presented in this paper are efficient and reliable, which is consistent with our theoretical analysis. But we have to note that the numerical eigenvalues do not perform  that well when $\mathbf{b}=(10,0)^T$. This is probably the
consequence of the performance of linear algebra routine on a convection dominated problem. \\

\indent {\bf Example 2.}~ Consider the convection-diffusion eigenvalue problem (\ref{s2.1}) on  $\Omega=(0,2)^2\setminus[1,2]^2$. Since the exact eigenvalues of (\ref{s2.1}) are unknown, we choose the approximate eigenvalues with high accuracy to replace them. For $\mathbf{b}=(1,0)^T, \mathbf{b}=(3,0)^T$, and $\mathbf{b}=(10,0)^T$, respectively, some adaptive refined meshes are shown in Figures 4 and 5 and the numerical results are shown in table 2. From the results we can see that for the convection parameters $\mathbf{b}=(1,0)^T, \mathbf{b}=(3,0)^T$, and
$\mathbf{b}=(10,0)^T$, the a posteriori error indicators can reflect the general trend of the error of discrete eigenvalues but similar to Example 1 the numerical eigenvalues do not  perform  that well when $\mathbf{b}=(10,0)^T$.

\begin{figure}
  \centering
   \scalebox{0.25}{\includegraphics{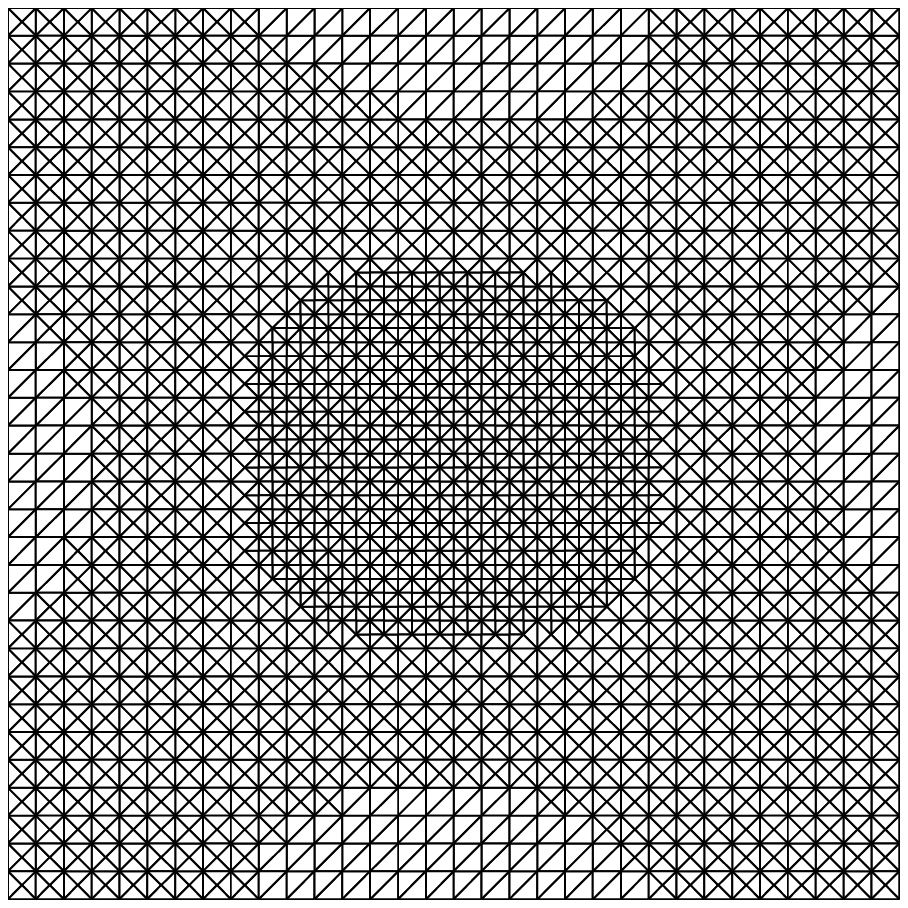}}
   \scalebox{0.25}{\includegraphics{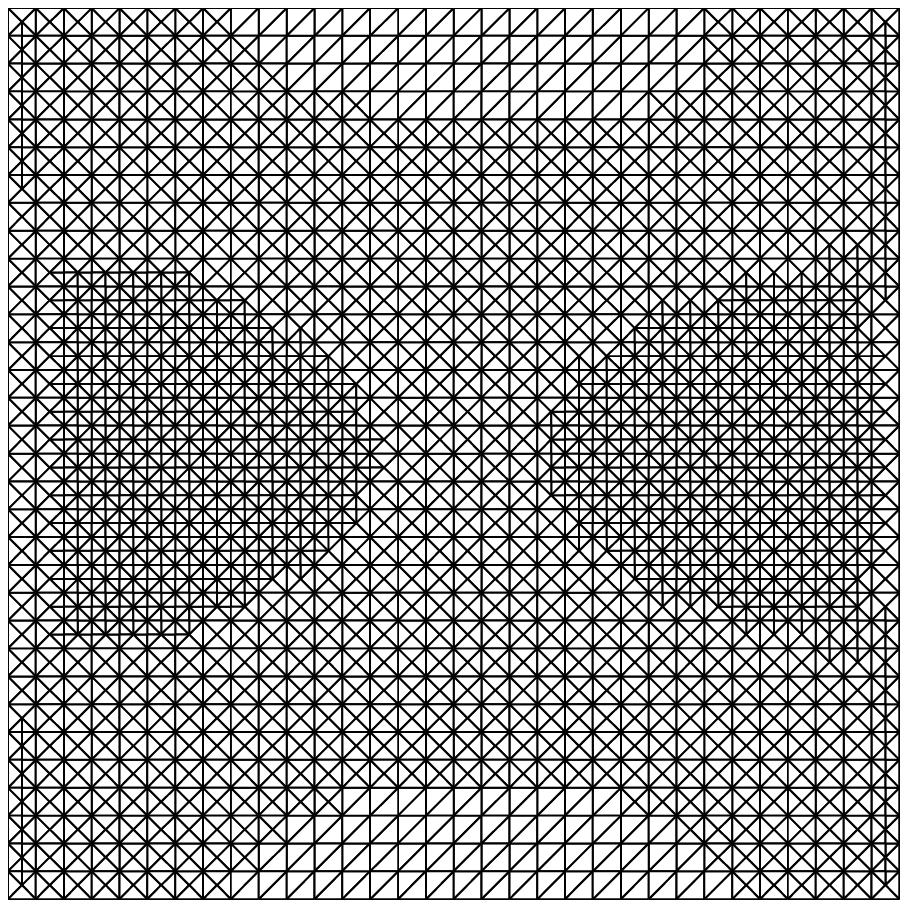}}
  \scalebox{0.25}{\includegraphics{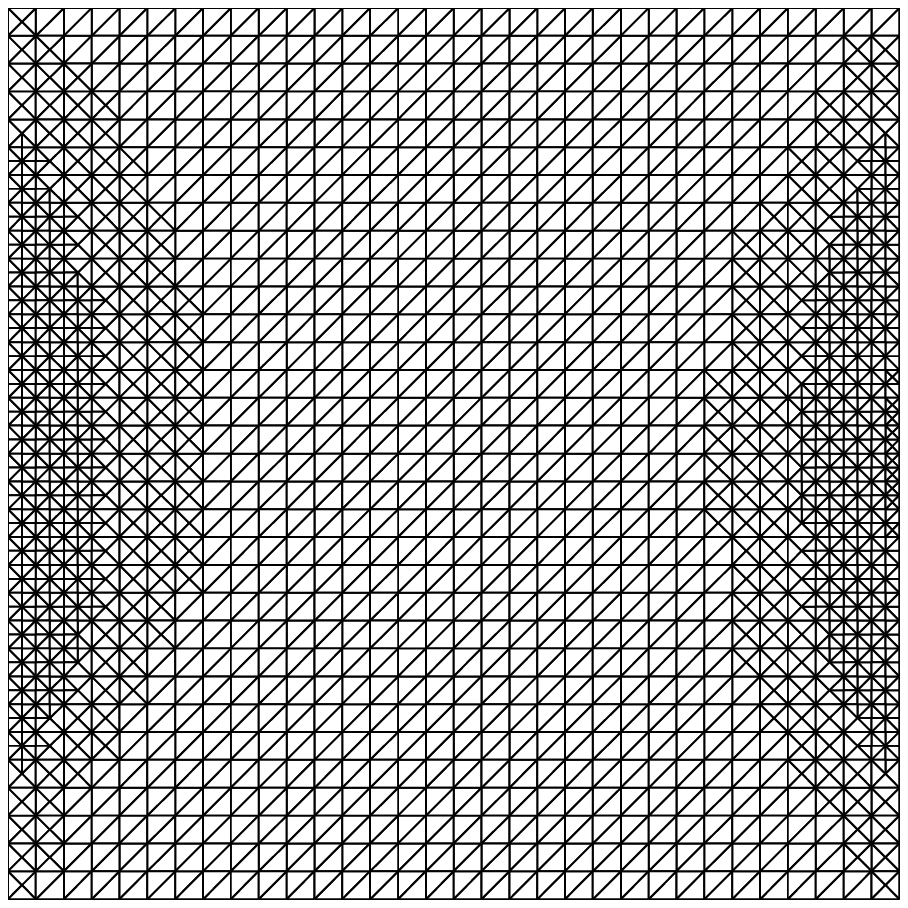}}
 \caption{the adaptively refined meshes of 1st eigenvalue after 6th iteration when $\mathbf{b}=(1,0)^T,\mathbf{b}=(3,0)^T$, and $ \mathbf{b}=(10,0)^T$, respectively.}
\end{figure}
\begin{figure}
  \centering
   \scalebox{0.25}{\includegraphics{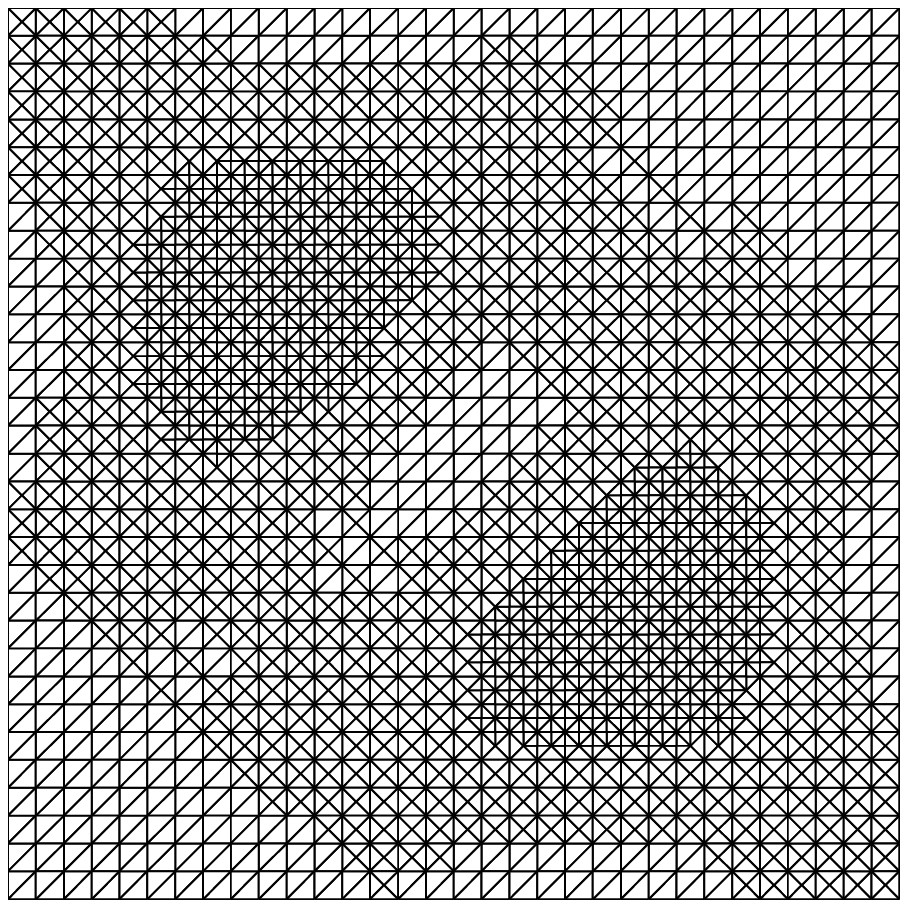}}
   \scalebox{0.25}{\includegraphics{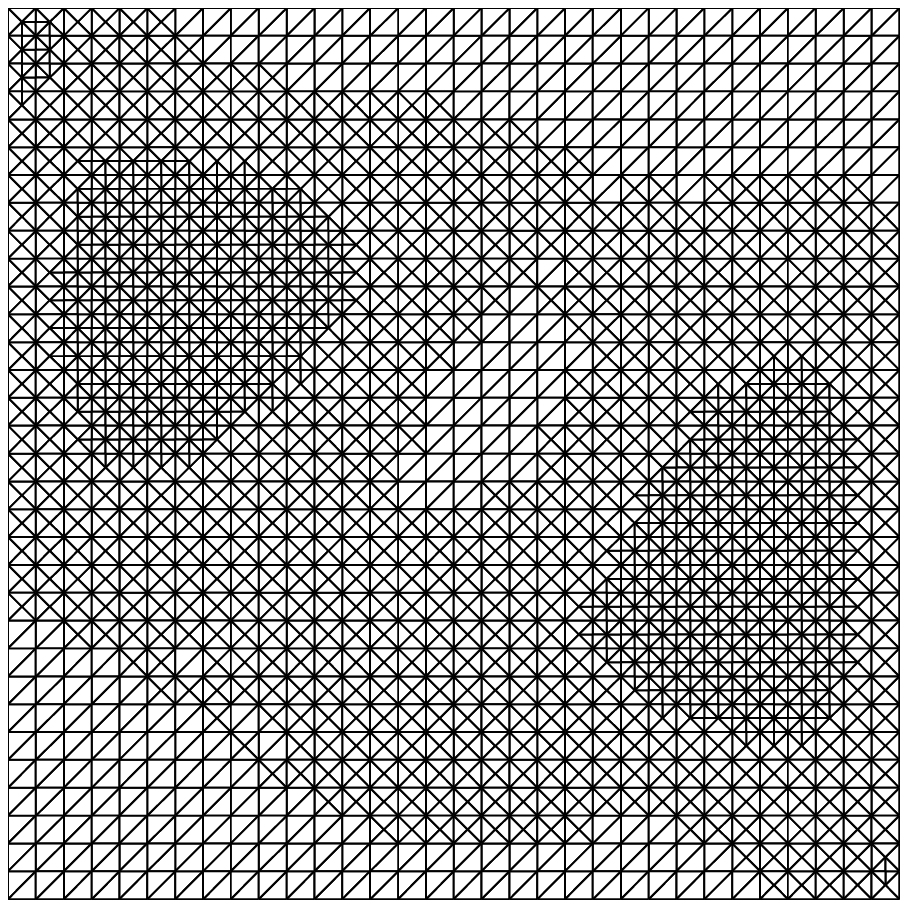}}
  \scalebox{0.25}{\includegraphics{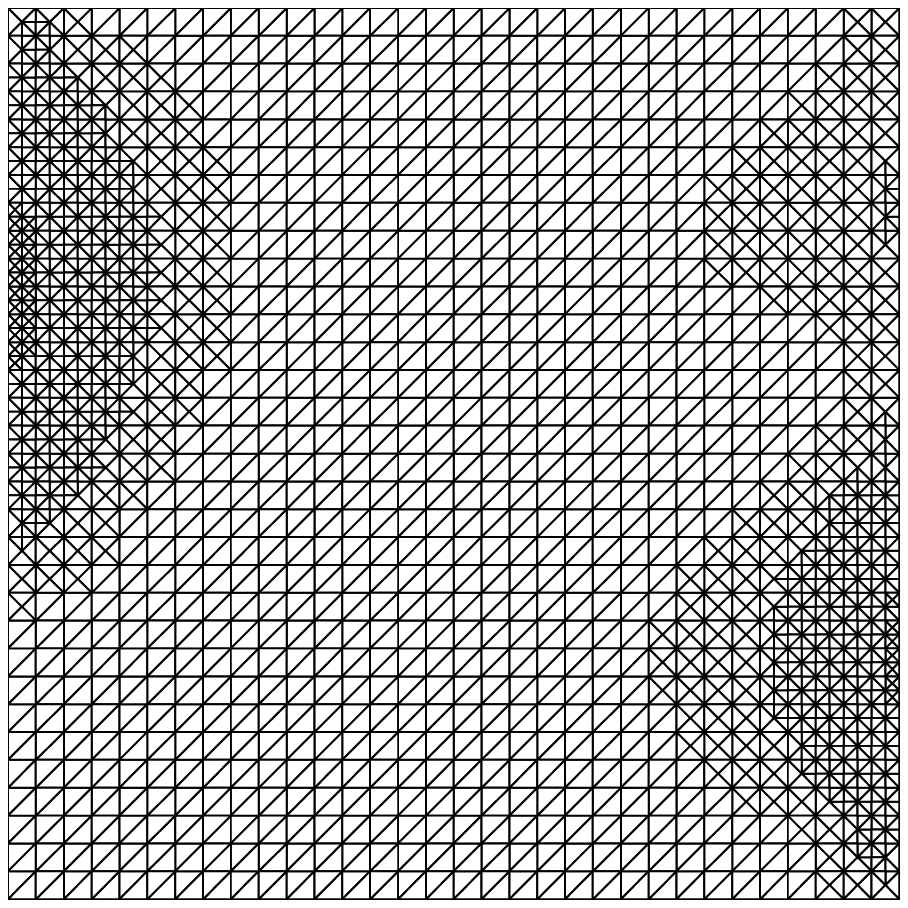}}
 \caption{the adaptively refined meshes of 2nd eigenvalue after 6th iteration when $\mathbf{b}=(1,0)^T, \mathbf{b}=(3,0)^T$, and $ \mathbf{b}=(10,0)^T$, respectively.}
\end{figure}
\begin{figure}
  \centering
   \scalebox{0.4}{\includegraphics{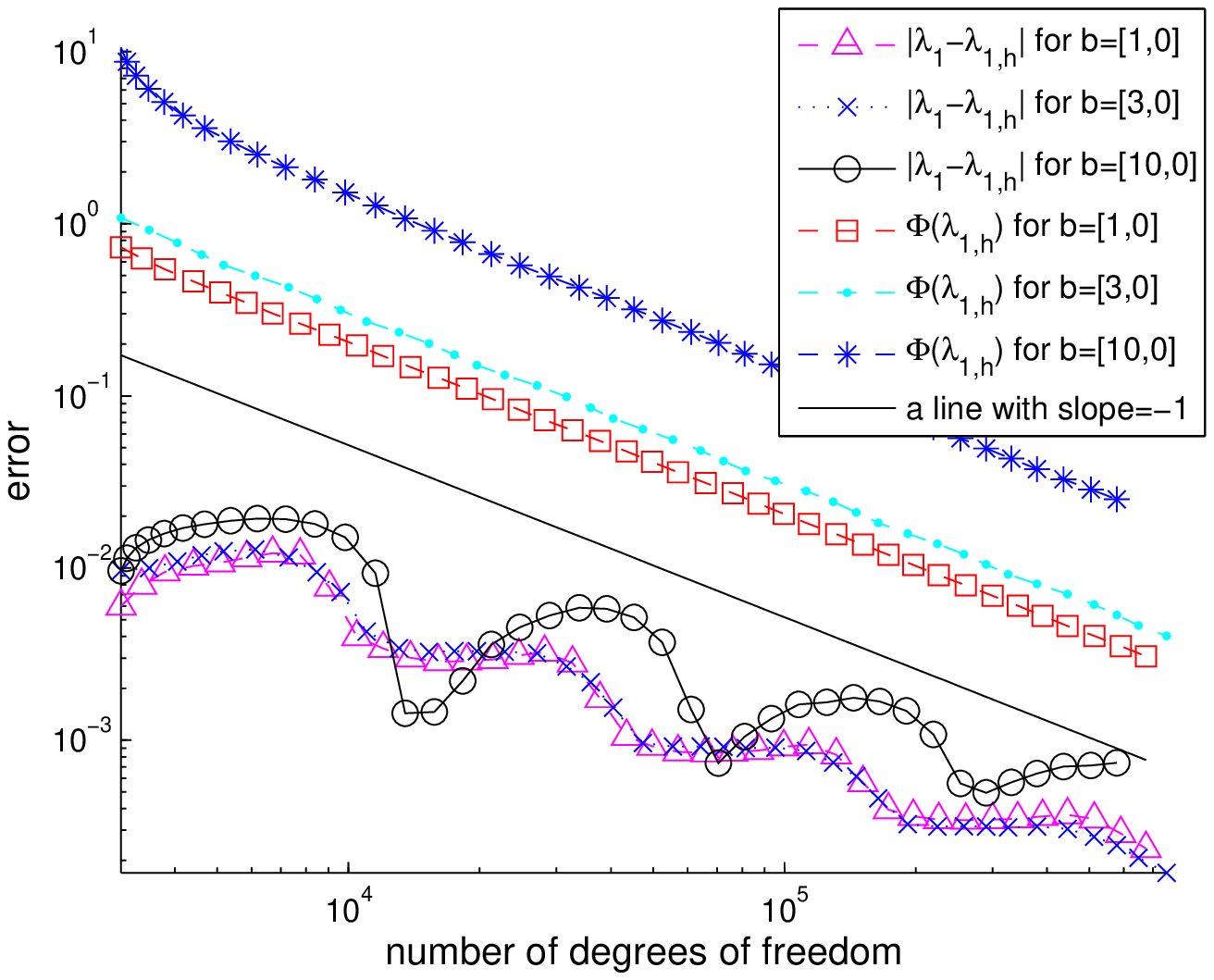}}
   \scalebox{0.4}{\includegraphics{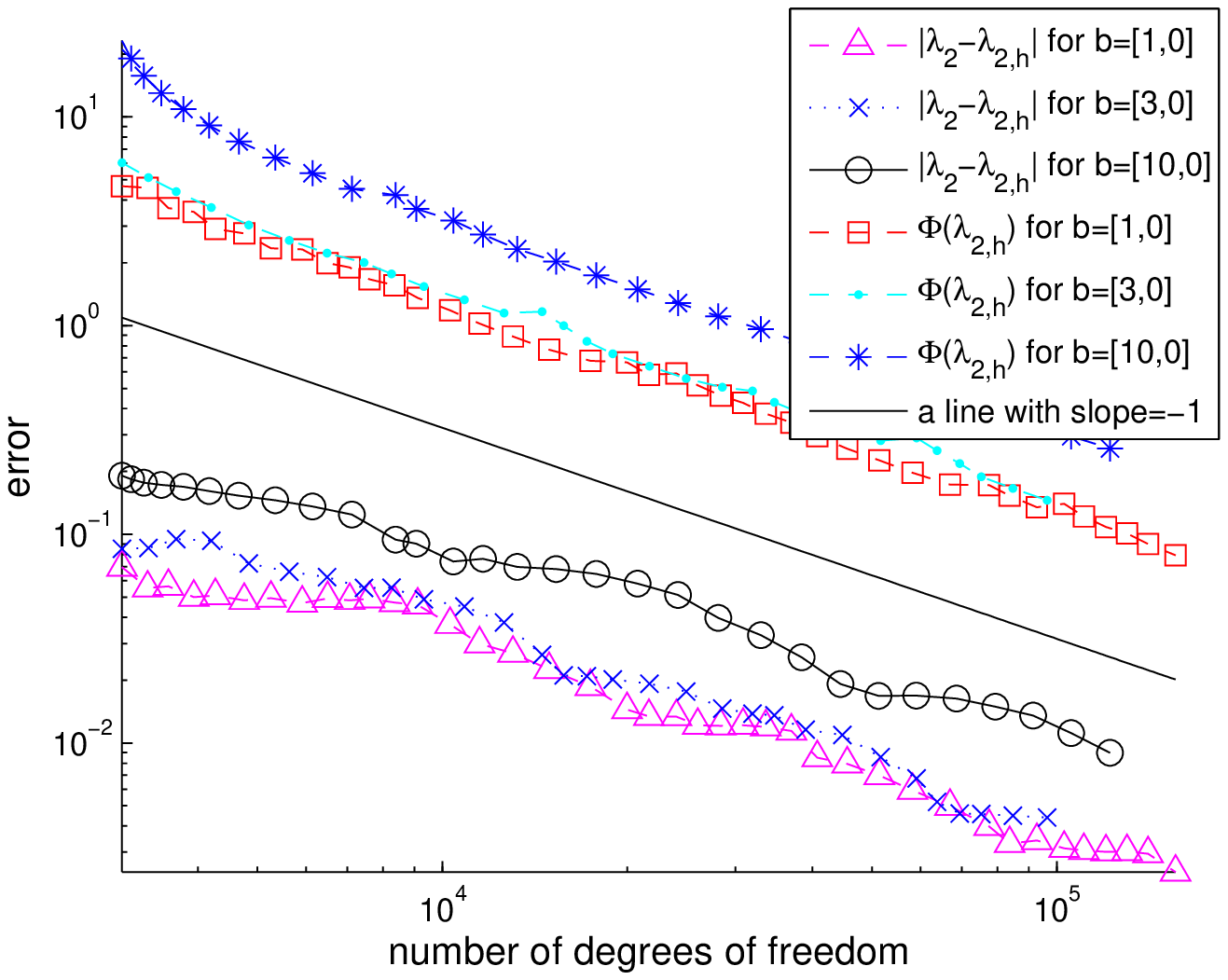}}
 \caption{$\Omega=(0,1)^2$, the first eigenvalue and the second eigenvalue}
\end{figure}
\begin{figure}
  \centering
   \scalebox{0.25}{\includegraphics{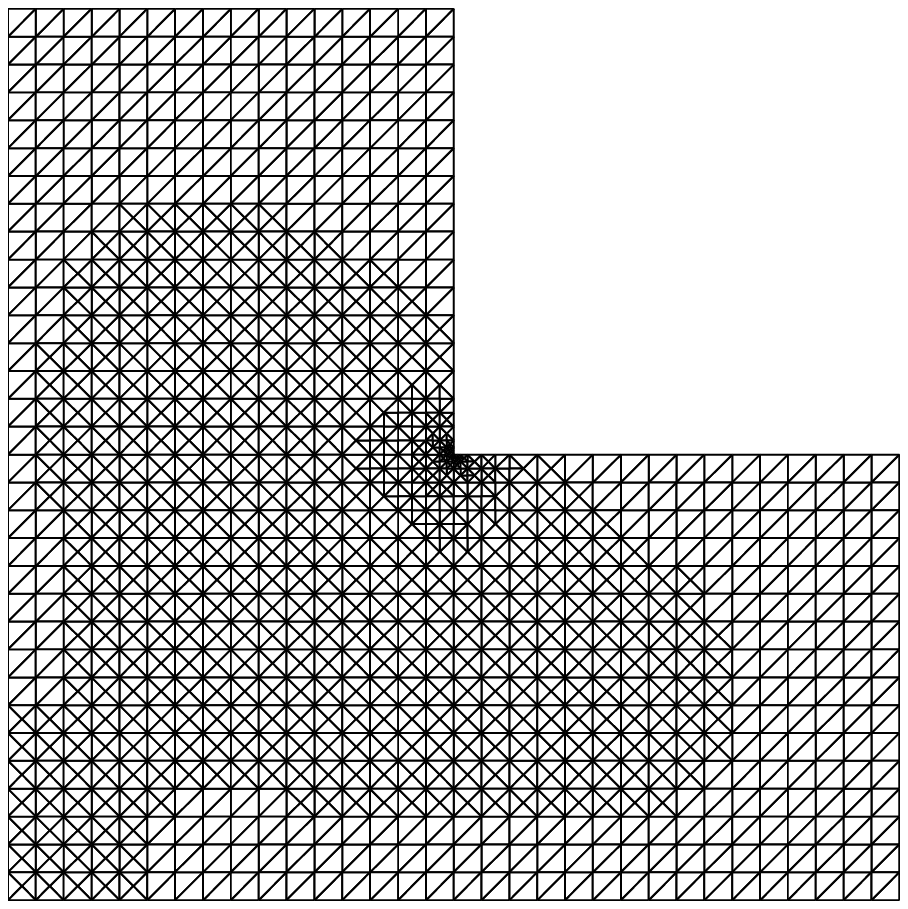}}
   \scalebox{0.25}{\includegraphics{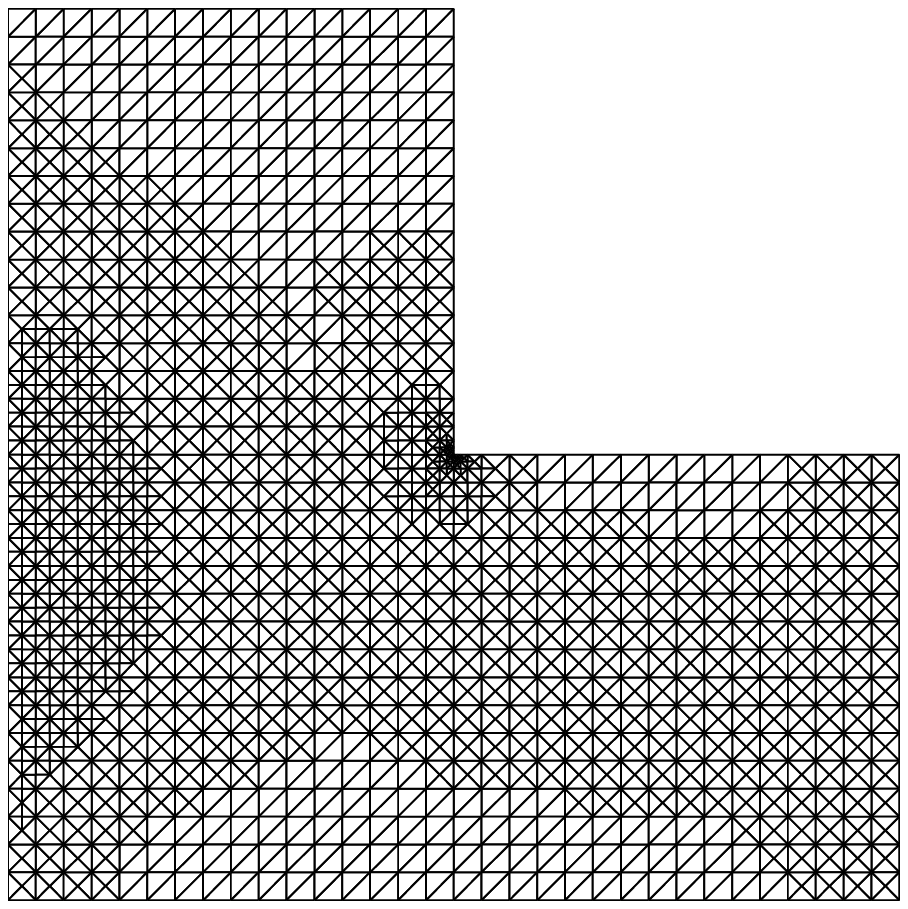}}
  \scalebox{0.25}{\includegraphics{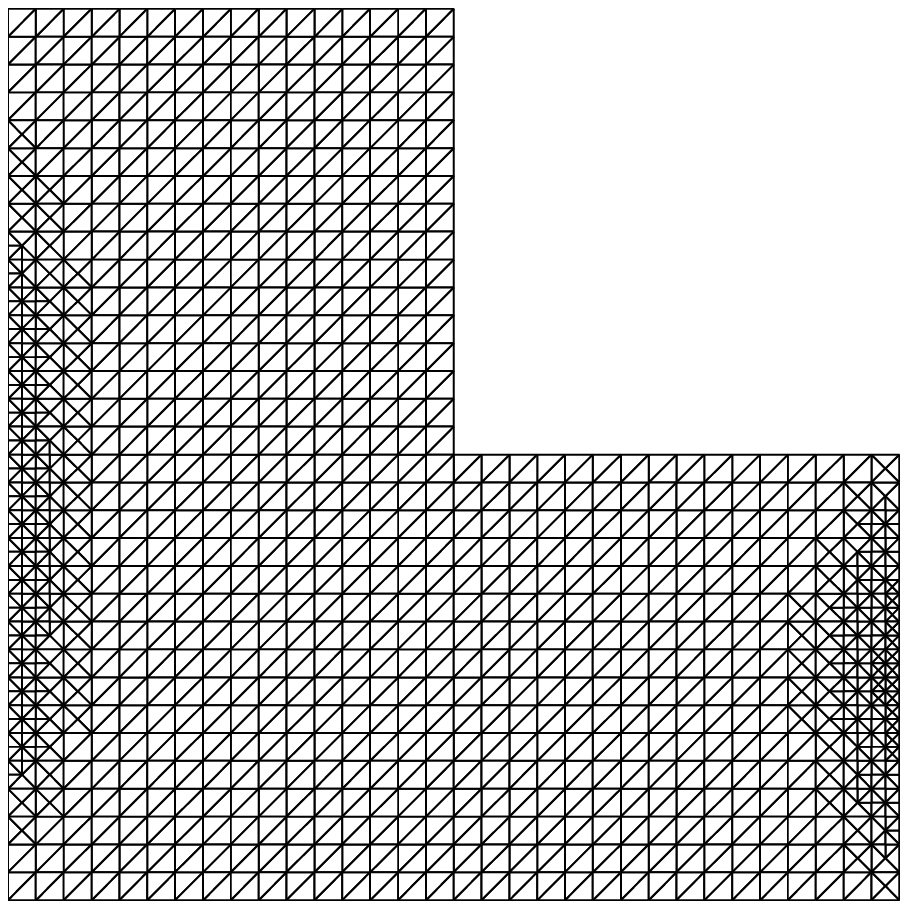}}
 \caption{the adaptively refined meshes of 1st eigenvalue after 6th iteration when $\mathbf{b}=(1,0)^T, \mathbf{b}=(3,0)^T$, and $ \mathbf{b}=(10,0)^T$, respectively.}
\end{figure}
\begin{figure}
  \centering
   \scalebox{0.25}{\includegraphics{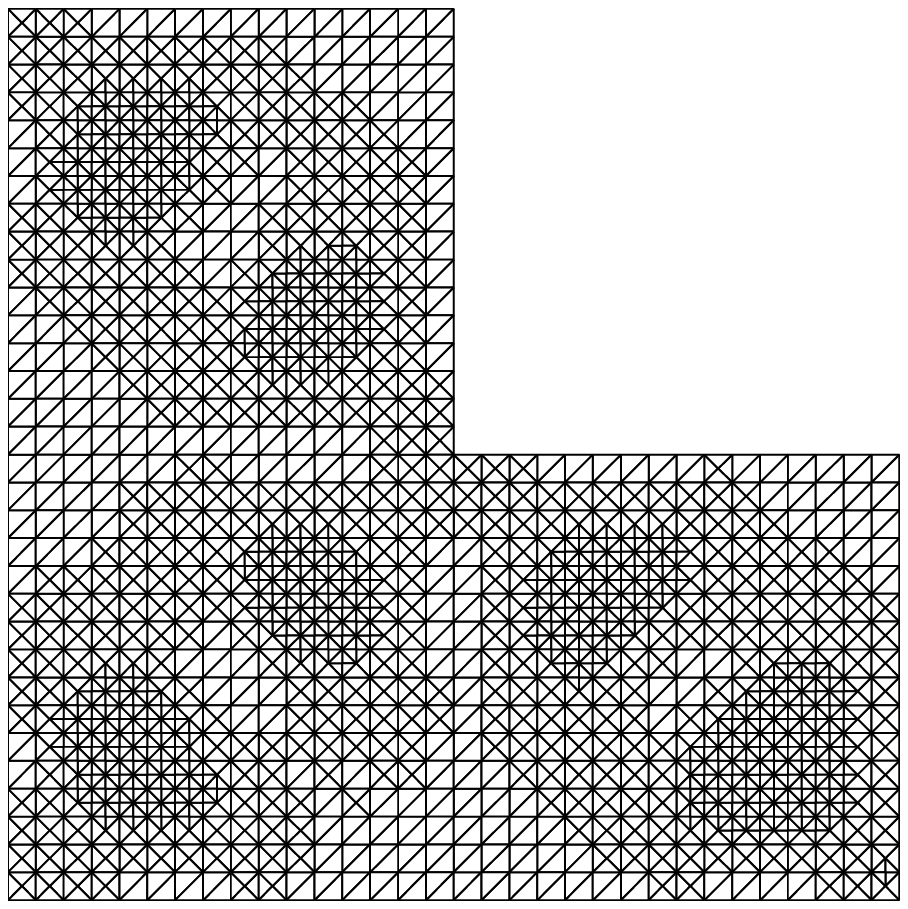}}
   \scalebox{0.25}{\includegraphics{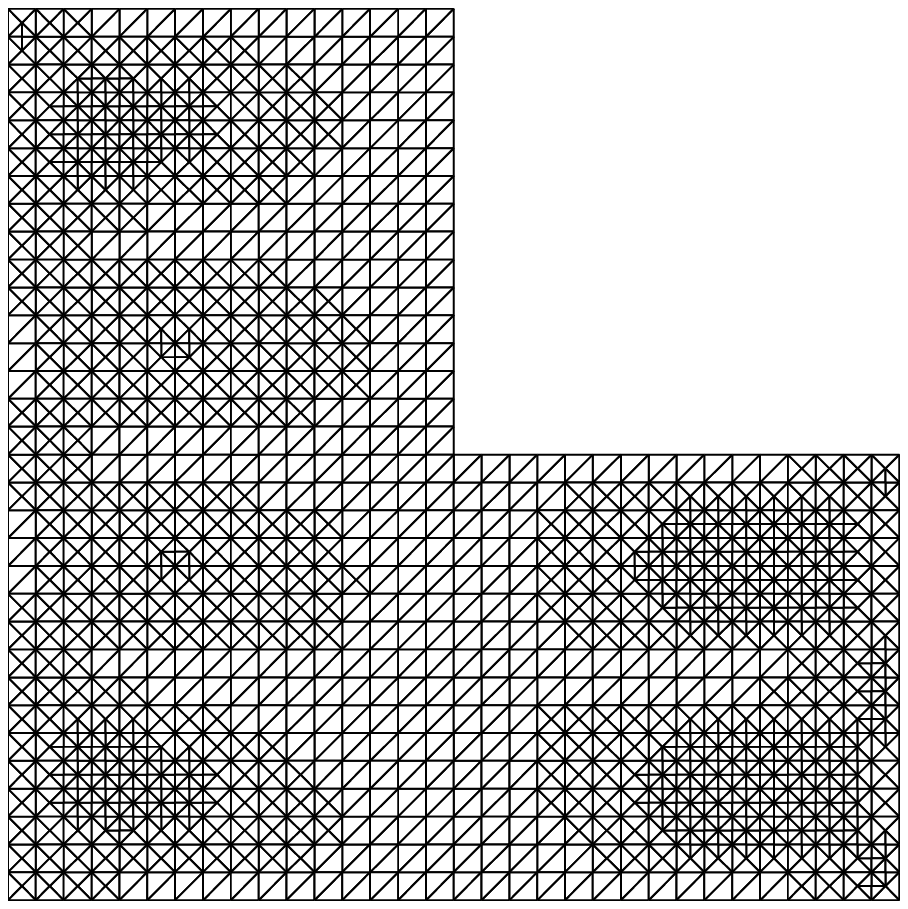}}
  \scalebox{0.25}{\includegraphics{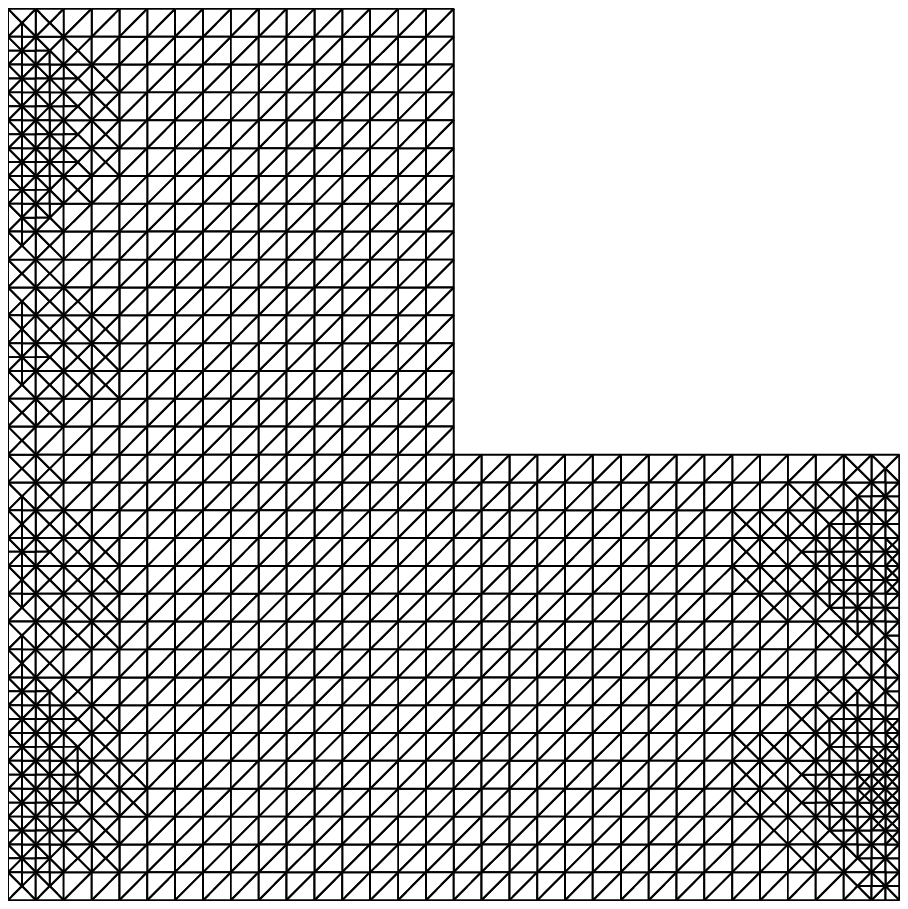}}
 \caption{the adaptively refined meshes of 8th eigenvalue after 6th iteration when $\mathbf{b}=(1,0)^T, \mathbf{b}=(3,0)^T$, and $ \mathbf{b}=(10,0)^T$, respectively.}
\end{figure}
\begin{figure}
  \centering
   \scalebox{0.4}{\includegraphics{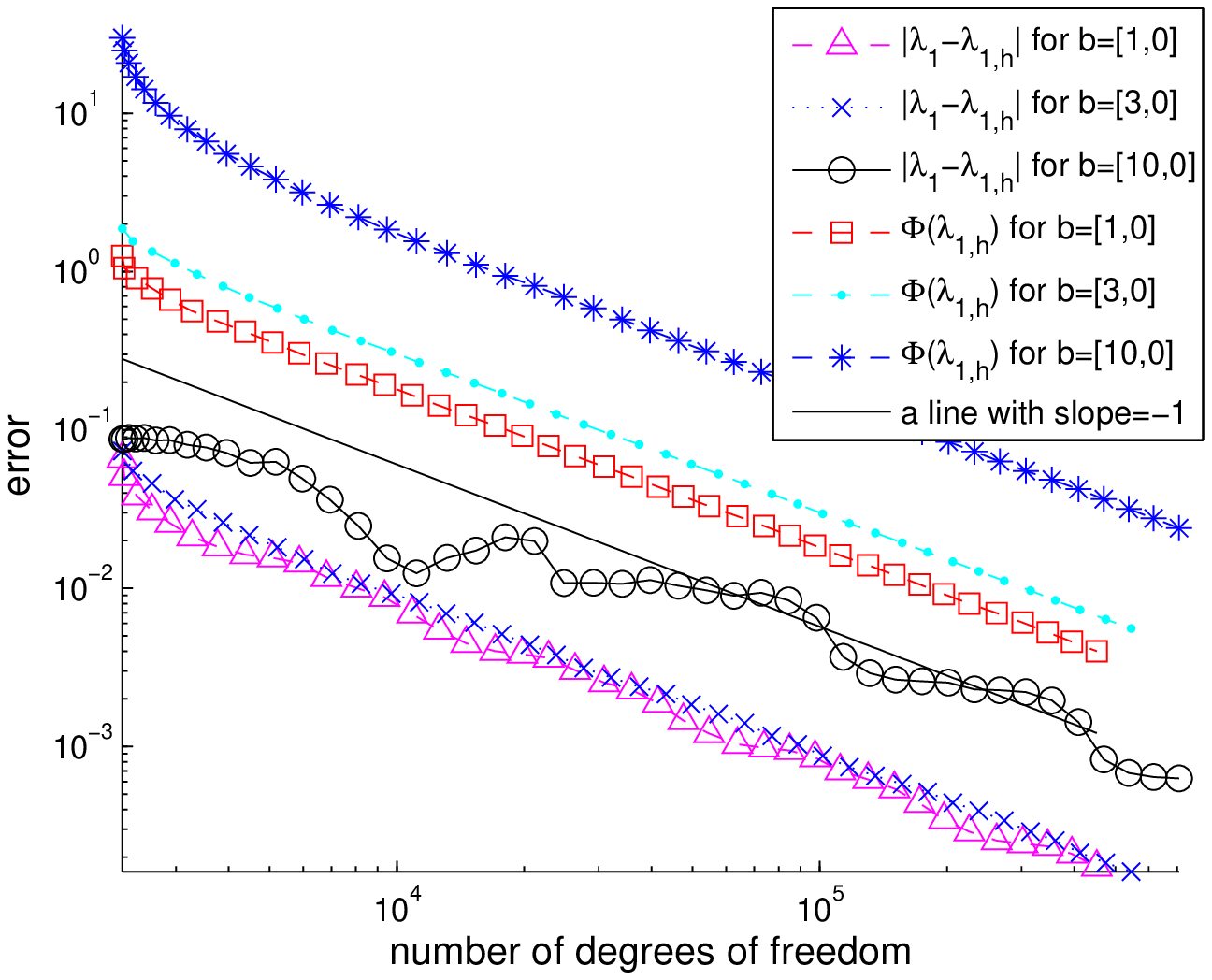}}
   \scalebox{0.4}{\includegraphics{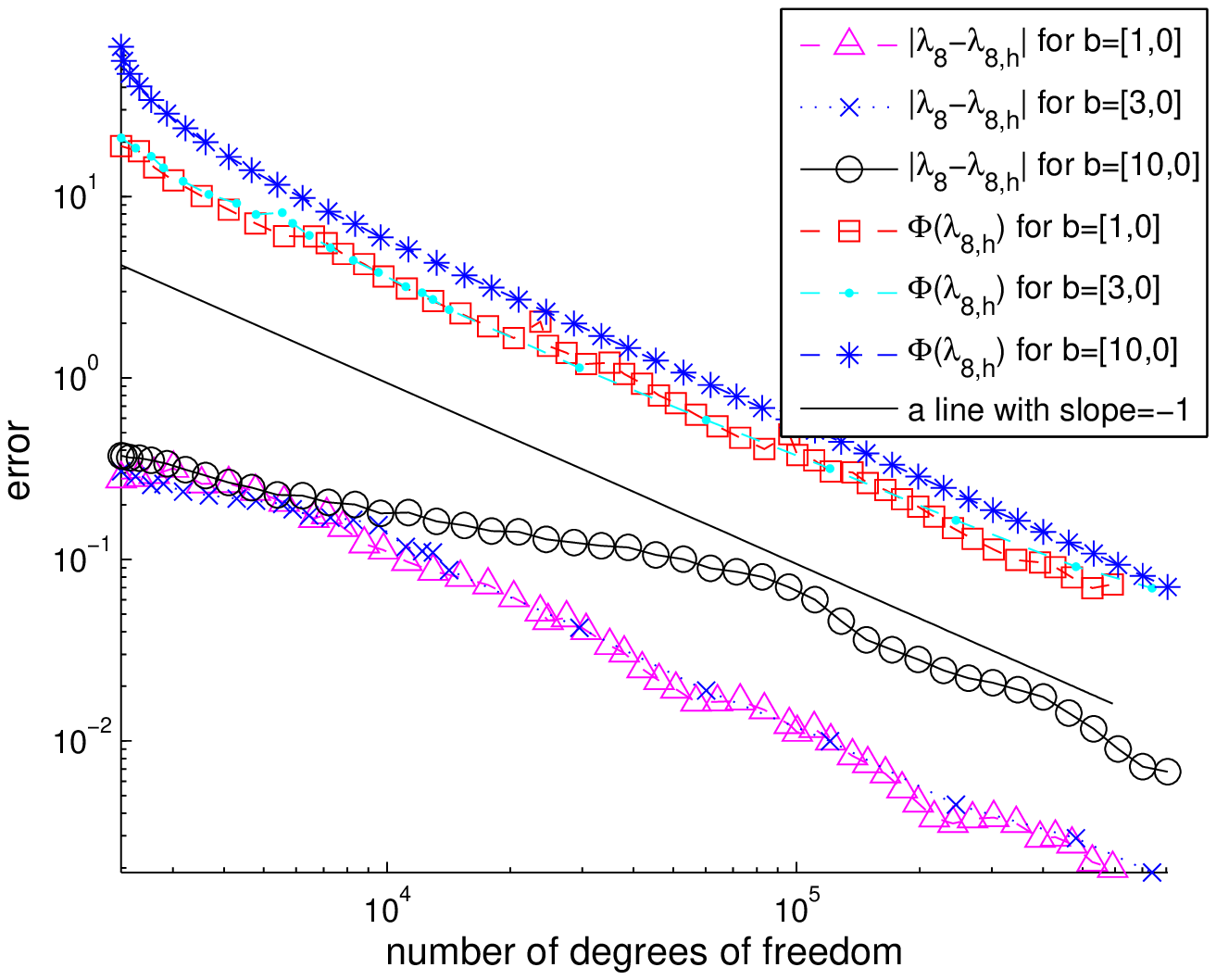}}
 \caption{$\Omega=(0,2)^2\setminus[1,2]^2$, the first eigenvalue and the 8th eigenvalue}
\end{figure}

\begin{table}
\caption{The 1st and 2nd eigenvalues on $\Omega=(0,1)^2$ with $H=\frac{\sqrt{2}}{16}$.}
\begin{center}\footnotesize
\begin{tabular}{ccccccccc}\hline
$k$&$l$&$N_{k,l}(1)$&$\lambda_{k}(1)$&$N_{k,l}(3)$&$\lambda_{k}(3)$&$N_{k,l}(10)$&$\lambda_{k}(10)$\\
\hline  1&6 &5846&19.977907&6113&21.976516&4171&44.755867\\
  1&18 &32624&19.986582&35943&21.987155&20781&44.743331\\
  1&30 &173175&19.988992&191642&21.989001&122039&44.740570\\
  1&38 &513308&19.989039&577955&21.989081&374065&44.739512\\
  1&39 &590647&19.989101&648352&21.989119&428658&44.739520\\
  1&40 &675033&19.989153&751651&21.989156&493913&44.739547\\
  2&6 &4757&49.54975734&5630&51.53186078&4174&74.18697239\\
  2&18 &17413&49.57934179&24916&51.58036988&24714&74.28997551\\
  2&30 &58202&49.59209659&96381&51.59363069&144114&74.33904514\\
  2&38 &130004&49.59500385&244514&51.59666467&436679&74.3446992\\
  2&39 &140739&49.59507564&266852&51.59682511&503989&74.34528594\\
  2&40 &155888&49.59561150&292042&51.59684387&578175&74.34579087\\
\hline
\end{tabular}
\end{center}
\end{table}
\begin{table}
\caption{The 1st and 8th eigenvalues on $\Omega=(0,2)^2\setminus[1,2]^2$ with $H=\frac{\sqrt{2}}{16}$.}
\begin{center}\footnotesize
\begin{tabular}{ccccccccc}\hline
$k$&$l$&$N_{k,l}(1)$&$\lambda_{k}(1)$&$N_{k,l}(3)$&$\lambda_{k}(3)$&$N_{k,l}(10)$&$\lambda_{k}(10)$\\
\hline  1&6 &3278&9.868593&3876&11.863678&2686&34.725519\\
  1&18 &19789&9.885896&23820&11.885939&13141&34.655247\\
  1&30 &111849&9.889023&135412&11.889074&84307&34.648104\\
  1&38 &345699&9.889490&412412&11.889511&266844&34.641987\\
  1&39 &394793&9.889514&474643&11.889540&307313&34.641930\\
  1&40 &452105&9.889551&545221&11.889562&354041&34.641676\\
  8&6 &4103&49.336511&3671&51.372524&2902&74.010829\\
  8&12 &8791&49.476412&7272&51.428181&6220&74.123357\\
  8&18 &20394&49.536447&14189&51.510658&15446&74.193993\\
  8&23 &34979&49.564431&481420&51.595113&33382&74.229178\\
  8&24 &37940&49.567477&738241&51.596140&38817&74.231469\\
  8&25 &41981&49.573008&1433565&51.596812&45354&74.242124\\
\hline
\end{tabular}
\end{center}
\end{table}

\begin{thebibliography}{s10}

\bibitem{J.G4}{J. Gedicke and C. Carstensen}, ``A posteriori error estimators for Convection-diffusion eigenvalue  problems,"
 Computer Methods in Applied Mechanics and Engineering, vol. 268, pp. 160-177, 2014.
\bibitem{R.R6}{R. Rannacher}, ``Adaptive FE eigenvalue computation with applications to hydrodynamic stability," in Advances in Mathematical Fluid Mechanics, pp. 425-450, Springer, Berlin, Germany, 2010.
\bibitem{YLHH7}{Y. Yang, L. Sun, H. Bi and H. Li}, ``A note on the residual type a posteriori error estimates for finite element eigenpairs of nonsymmetric elliptic eigenvalue problems," Applied Numerical Mathematics, vol. 82, pp. 51-67, 2014.
\bibitem{C.Carstensen4}{C. Carstensen, J. Gedicke, V. Mehrmann and A. Miedlar}, ``An adaptive homotopy approach for non-selfadjoint eigenvalue problems," Numerische Mathematik, vol. 119, no. 3, pp. 557-583, 2011.
\bibitem{TLYF9}{$T. L\ddot{u}$ and Y. Feng}, ``Splitting extrapolation based on domain decomposition for finite element approximations," Science in China Series E:Technological Sciences, vol. 40, no. 2, pp. 144-155, 1997.
\bibitem{YHS10}{Y. Yang, H. Bi and S. Li}, ``The extrapolation of numerical eigenvalues by finite elements for differential operators," Applied Numerical Mathematics, vol. 69, pp. 59-72, 2013.

\bibitem{A.N}{A. Naga and Z. Zhang}, ``Function value recovery and its application in eigenvalue problems,"
SIAM Journal on Numerical analysis, vol. 50, no. 1, pp. 272-286,
2012.

\bibitem{JHYY5}{J. Han and Y. Yang}, ``A class of spectral element methods and its a priori/a posteriori error estimates
 for 2nd-order elliptic eigenvalue problems," Abstract and Applied Analysis, vol. 2013, 14 pages, 2013.

\bibitem{YH}{Y. Yang and J. Han}, ``Multilevel finite element discretizations based on local defect
correction for nonsymmetric eigenvalue problems," Computers and
Mathematics with Applications, vol. 70, pp. 1799-1816, 2015.

\bibitem{pzlan}{Z. Peng, H. Bi, H. Li and Y. Yang}, ``A Multilevel correction method for Convection-diffusion eigenvalue  problems,"
Mathematical Problems in Engineering, Vol. 2015, pp. 1-10, 2015.

\bibitem{ainsworth1}{M. Ainsworth and J.T. Oden},``A posterior error estimation in Finite element Analysis,"
 Wiley-Inter science, New York, 2011.

\bibitem{babuska3}{I. Babuska and W.C. Rheinboldt}, ``Error estimates for adaptive finite element
computations," SIAM J.Numer.Anal., Vol. 15, pp. 736-754, 1978.

\bibitem{verfurth1} {R. Verf$\ddot{u}$rth}, ``A Posteriori Error Estimation Techniques,"
Oxford University Press, USA, 2013.

\bibitem{shi2}{Z. Shi and M. Wang}, ``Finite Element Methods," Beijing, Scientific Publishers, 2013.

\bibitem{C. Carstensen2}{C. Carstensen, J. Hu and A. Orlando}, ``Framework for the a posteriori error analysis of noncomforming finite elements,"
SIAM. J. Numer. Anal, vol. 45, pp. 68-82, 2007.

\bibitem{Li}{Y. Li}, ``A posteriori error analysis of nonconforming methods for the eigenvalue problem," Jrl Syst Sci \& Complexity,
 vol. 22, pp. 495-502, 2009.

\bibitem{YHB}{Y. Yang, J. Han and H. Bi}, ``Non-conforming finite element methods for transmission eigenvalue problem," Computer Methods in Applied Mechanics and Engineering, vol. 307, pp. 144-163, 2016.

\bibitem{C. Carstensen8}{C. Carstensen and J. Hu}, ``A unifying theory of a posteriori error control for nonconforming finite element methods," Numer. Math, vol. 107, pp. 473-502, 2007.
\bibitem{R.V}{R. Verf$\ddot{u}$rth}, ``A Review of a Posteriori Error Estimation and Adaptive Mesh-Refinement Techniques" Wiley-Teubner, 2007.

\bibitem{P.C}{$\mathrm{P. Cl\acute{e}ment}$}, ``Approximation by finite functions using local regularization," RAIRO Ana Numer, vol. 9, pp. 77-84, 1975.

\bibitem{C.Carstensen3}{C. Carstensen}, ``Quasi-interpolation and a poeteriori error analsis in finite element methods,"
 ${M}^2$NA, vol. 33, pp. 1187-1202, 1999.

\bibitem{C.B}{C. Bernardi and V. Girault}, ``A local regularisation operator for triangular and quadrilateral finite elements,"
SIAM J. Numer. Anal, vol. 35, pp. 1893-1916, 1998.

\bibitem{L.chen}{L. chen}, ``iFEM: an innovative finite element methods package in
MATLAB," Technical Report, University of California at Irvine, 2009.


\end{thebibliography}
\end{document}